\documentclass[11pt]
{article}
\usepackage{amsmath,amsfonts,amssymb,amsthm} 
\usepackage {amsmath, amssymb,amsthm, graphicx, mathrsfs}
\usepackage{amsmath, amssymb, amsthm, authblk, bigints}
\theoremstyle {plain}
\usepackage{color}
\usepackage{enumerate}

\newtheorem{lemma}{Lemma}

\newtheorem{theorem}{Theorem}
\newtheorem{remark}{Remark}
\newtheorem{proposition}{Proposition}

\newtheorem{definition}{Definition}
\newtheorem{assumption}{Assumption}

\newcommand{\ho}{\mathrm{hom}}
\newcommand{\dist}{\mathrm{dist}}
\newcommand{\SO}{\mathrm{SO}}
\newcommand{\R}{\mathbb{R}}
\newcommand{\Z}{\mathbb{Z}}

\newcommand{\e}{\varepsilon}
\mathchardef\emptyset="001F
\newcommand{\wto}{\rightharpoonup} 
\newcommand{\wtto}{\stackrel{2}{\rightharpoonup}} 
\newcommand{\stto}{\stackrel{2}{\rightarrow}}
\newcommand{\sym}{{\operatorname{sym}}}

\newcommand{\step}[1]{\noindent \textit{Step} #1.}

\begin{document} 
\title{\sc High contrast homogenisation in nonlinear elasticity under small loads}
\author[1]{Mikhail Cherdantsev}
\author[2]{Kirill Cherednichenko} 
\author[3]{Stefan Neukamm} 
\affil[1]{School of Mathematics, Cardiff University, Senghennydd Road, Cardiff, CF24 4AG, UK}
\affil[2]{Department of Mathematical Sciences, University of Bath, Claverton Down, Bath, BA2 7AY, 
UK} 
\affil[3]{TU Dresden, Fachrichtung Mathematik, 
Institut f\"{u}r Wissenschaftliches, Rechnen, D-01062 Dresden, Germany}

\maketitle

\begin{abstract}
We  study the homogenisation of \textit{geometrically nonlinear elastic}
composites with high contrast. The composites we analyse consist of a
perforated matrix material, which we call the ``stiff'' material, and
a ``soft'' material that fills the pores. We assume that the pores are
of size $0<\e\ll 1$ and are periodically distributed with period $\e$.
We also assume that the stiffness of the soft material
degenerates with rate $\e^{2\gamma},$ $\gamma>0$, so that
the contrast between the two materials becomes
infinite as $\e\downarrow 0$. We study the homogenisation limit
$\e\downarrow 0$ in a low energy regime, where the displacement of the stiff component is infinitesimally small. We derive 
an effective
two-scale model, which, depending on the scaling of the energy, is
either a quadratic functional or a partially quadratic functional that 
still allows for large
strains in the soft inclusions. 
In the latter case, averaging out the small scale-term justifies a single-scale model for high-contrast materials, which features a non-linear and non-monotone effect describing a coupling between microscopic and the effective macroscopic displacements.
\medskip

\end{abstract}

\tableofcontents

\section{Introduction}
\label{intro_section}

We consider a geometrically nonlinear elastic composite material that consists of a ``stiff'' matrix material and periodically distributed pores filled by a ``soft'' material: for $\e>0$ and a fixed scaling parameter $\gamma>0$ we consider the energy functional of non-linear elasticity
\begin{equation}\label{eq:energy:00}
  \mathcal I_\e(u):=\int_\Omega
  \Bigl(\e^{2\gamma}W^0(\nabla u)\chi_\e+W^1(\nabla u)(1-\chi_\e)\Bigr)\,dx-\int_\Omega f_\e\cdot u\,dx,\qquad
  u\in H^1(\Omega).
\end{equation}
Here $\Omega$ denotes a Lipschitz domain in $\R^d$ (the
reference domain of the elastic body) and $u:\Omega\to\R^d$ is a deformation satisfying clamped boundary conditions: $u(x)=x$ on $\partial\Omega$. We denote by $f_\e:\Omega\to\R^d$ the density of the applied body forces, $W^0$ and $W^1$ are
frame-indifferent, non-degenerate energy densities (see Section~\ref{S:ass} below for
the precise assumptions), and $\chi_\e$ denotes the indicator function
of the pores, {\it i.e.}~of the domain occupied by the ``soft'' material component. As will
be made precise in Section~\ref{S:ass}, we assume that the pores are of size $\e$ and are
periodically distributed in the interior of $\Omega$ with period $\e$.
 As can be seen from (\ref{eq:energy:00}), in the homogenisation limit $\e\downarrow 0$ the stiffness of the ``soft'' material degenerates
with rate $\e^{2\gamma}$ $(\gamma>0$), while the stiffness of the ``stiff'' material
remains unchanged. Hence, the contrast between the soft
material (occupying the pores) and the stiff material (occupying
the perforated matrix) becomes infinite in the limit $\e\downarrow 0.$ We therefore refer to the corresponding limit procedure as \textit{high-contrast homogenisation.} Our goal is to identify the effective behaviour of $\mathcal I_\e$ by studying its limit under a proper rescaling. 

\paragraph{Summary and discussion of our result.} To illustrate our result, here in the introduction we restrict ourselves to the special case $\gamma=1$. If we assume that the density of the body forces is small in magnitude, in the sense that $f_\e=\e^{\alpha} f$ for some $\alpha \geq 1$, and has vanishing first moment, {\it i.e.} $\int_{\Omega} f(x)\cdot x\,dx=0$, then \eqref{eq:energy:00} can be expressed as
\begin{equation*}
  I_\e^\alpha(\varphi):=\frac{1}{\e^{2\alpha}}\mathcal I_\e(u)=\int_\Omega
  \biggl(\frac{1}{\e^{2(\alpha-1)}}W^0(I+\e^\alpha\nabla \varphi)\chi_\e+\frac{1}{\e^{2\alpha}}W^1(I+\e^{\alpha}\nabla \varphi)(1-\chi_\e)\biggr)\,dx-\int_\Omega f\cdot \varphi\,dx,
\end{equation*}
where 
\begin{equation}
\varphi(x)=\frac{u(x)-x}{\e^\alpha},\ \ \ x\in\Omega,\ \ \ \ \ \ \ \varphi\in H^1_0(\Omega;\R^d),
\label{u_to_phi}
\end{equation}
denotes the scaled displacement, and $I$ stands for the identity matrix in $\R^{d\times d}$. In this paper we prove (in fact, in a more general situation) that the functionals $I_\e^\alpha$ two-scale $\Gamma$-converge as $\e\downarrow 0$. In doing so, we distinguish two regimes, as follows.
\medskip

In the \textbf{small strain regime},  which corresponds to $\alpha>1$, the strain $\e^\alpha\nabla\varphi$ becomes infinitesimally small in the
  entire domain $\Omega$, and the limit behaviour is expressed by a linearised, two-scale
  energy acting on pairs of functions $(g^0, g^1)$, where $g^1\in H^1_0(\Omega;\R^d)$
  describes the (scaled) macroscopic displacement of the body, and
  $g^0\in L^2\bigl(\Omega;H^1_0(Y_0)\bigr)$ is a two-scale function describing the (scaled) microscopic displacement on the pores
  relative to the deformed matrix material. The two-scale $\Gamma$-limit is given by the convex quadratic functional
  \begin{eqnarray}
    I_{\textrm{small}}(g^0,g^1)&:=&\int_{\Omega\times Y}Q^0\bigl(\nabla_y g^0(x,y)\bigr)\,dy+Q^1_{\hom}\bigl(\nabla g^1(x)\bigr)\,dx\label{first_term}\\
    &&-\int_\Omega\left(\int_{Y_0}g^0(x,y)\,dy+g^1(x)\right)\cdot f(x)\,dx.\label{second_term}
  \end{eqnarray}
  Here $Q^0$ and $Q^1$ are the quadratic forms of the quadratic expansions of $W^0$ and $W^1$ at the identity, and $Q^1_{\hom}$ denotes the homogenised energy density obtained from $Q^1$, see \eqref{WTaylor} and \eqref{def:Q1hom} for details.
  
  We further prove ({\it cf.} Theorem~\ref{T1b}) that if $\varphi^\e$ is an (almost) minimiser of $I^0_\e$, then it admits the two-scale expansion
  \begin{equation}\label{intro:small:exp}
    \varphi_\e(x)\approx g^1_*(x)+\e\psi(x,x/\e)+g^0_*(x,x/\e),
  \end{equation}
  where $\psi$ denotes a corrector function that only depends on $g^1_*,$ and $(g^0_*,g^1_*)$ is a minimiser of $I_{\text{small}}$.

  By averaging out the fast variable $y$, the limit $I_{\textrm{small}}$ can be further simplified. In fact, we show that $I^\alpha_\e$ $\Gamma$-converge to the functional $\bar I_{\text{small}}:L^2(\Omega)\to \R\cup\{+\infty\}$ given by
  \begin{eqnarray}\label{intro:small:weak-single-scale}
    \bar I_\text{small}(\varphi)&:=&\min\left\{\int_\Omega Q^1_{\hom}(\nabla g^1)+\overline Q(\overline G)\,:\,g^1\in H^1_0(\Omega),\,\overline G\in L^2(\R^d)\text{ with }g^1+\overline G=\varphi\,\right\}\\\notag
                                 &&-\int_\Omega f\cdot\varphi\,dx,
  \end{eqnarray}
  with a quadratic form $\overline Q:\R^d\to[0,\infty)$, which is defined by
  \begin{equation*}
    \overline Q(\overline G):=\min\left\{\int_{Y^0}Q^0\bigl(\nabla g^0(y)\bigr)\,dy\,:\,g^0\in H^1_0(Y^0),\,\int_{Y^0}g^0\,dy=\overline G\right\},
  \end{equation*}
  and which captures the influence of the pores (and their geometry) on the effective behavior. The minimiser $\varphi_*\in L^2(\Omega)$ to $\bar I^\alpha_0$  takes the form $\varphi_*= g^1_*+\overline G_*$ with $\overline G_*:=\int_{Y_0}g^0_*(\cdot,y)\,dy$.
  In view of \eqref{intro:small:exp}  the field $\overline G_*$ can be interpreted as the gap between the macroscopic displacement and the microscopic displacements in the pores.  See  Remark~ \ref{R:intro:small} for more details.
\medskip

\medskip

In the  \textbf{finite strain regime}, which corresponds to $\alpha=1$, the displacement gradient $\e\nabla\varphi_\e$ becomes infinitesimally small only in the stiff component, while large strains still may occur in the soft pores. Therefore, the $\Gamma$-limit is a non-convex (partially linearised) functional of the form
\begin{eqnarray*}
  I_{\text{\rm finite}}(g^0,g^1)&:=&\int_{\Omega\times Y}\mathcal QW^0\bigl(I+\nabla_y g^0(x,y)\bigr)\,dy+Q^1_{\hom}\bigl(\nabla g^1(x)\bigr)\,dx\\
  &&-\int_\Omega\left(\int_{Y_0}g^0(x,y)\,dy+g^1(x)\right)\cdot f(x)\,dx,
\end{eqnarray*}
where $\mathcal QW^0$ denotes the quasiconvex envelope of $W^0$. Similarly to the small strain regime, one can 
average out the fast scale $y$ and obtain $\Gamma$-convergence of $I^0_\e$ to the functional $\bar I_{\text{finite}}:L^2(\Omega)\to\R\cup\{+\infty\}$ given by
\begin{eqnarray}\label{intro:finite:weak-single-scale}
  \bar I_\text{\rm finite}(\varphi)&:=&\min\left\{\int_\Omega Q^1_{\hom}(\nabla g^1)+\overline V(\overline G)\,:\,g^1\in H^1_0(\Omega),\,\overline G\in L^2(\R^d)\text{ with }g^1+\overline G=\varphi\,\right\}\\\notag
                              &&-\int_\Omega f\cdot\varphi\,dx,
\end{eqnarray}
with non-convex potential $\overline V:\R^d\to[0,\infty)$ defined by
\begin{equation*}
  \overline V(\overline G):=\min\left\{\int_{Y^0}\mathcal QW^0\bigl(I+\nabla g^0(y)\bigr)\,dy\,:\,g^0\in H^1_0(Y^0),\,\int_{Y^0}g^0\,dy=\overline G\right\},
\end{equation*}
see  Remark~\ref{R:intro:finite}. In contrast to the small strain regime, where $\overline Q$ is quadratic, the potential $\overline V$ is non-convex and expresses a \textit{nonlinear (and non-monotone) coupling} between the macroscopic and microscopic displacement.
\medskip

\paragraph{Connection to acoustic wave propagation in high-contrast materials.} The sequence of functionals $\varepsilon^{-2\alpha}{\mathcal I}_\varepsilon,$ in either of the two regimes described above, occupies an intermediate position between a fully nonlinearly elastic composite and fully linearised models, as $\varepsilon\to0.$ Notably, linear models with high contrast, which are suitable for the description of small displacement fields (that often occur, say, in acoustic wave propagation) already exhibit a coupling between the macroscopic part $g^0$ and microscopic part $g^1$ of the minimiser of $I_{\rm small},$ which in our case is obtained as a limit in the small-strain regime $\alpha>1.$ 
This can be seen by considering the time-harmonic solitons to the equations of elastodynamics with the elastic part of the energy given by (\ref{first_term}), away from the sources of the elastic motion. In this case 
the function $f$ in (\ref{second_term}) (which in our analysis we assume to be independent of the fast variable $x/\varepsilon$ for simplicity, an assumption that can be relaxed with no changes in the proofs needed) has to be replaced by the sum $g^0+g^1,$ with the integration in (\ref{second_term}) carried over 
$Y_0$ and $Q$ at the same time, {\it i.e.} the work of the external forces (\ref{second_term}) is replaced by the expression for the work of ``self-forces'' 
$$
-\omega^2\int_\Omega\int_{Y_0}\left(g^0(x,y)+g^1(x)\right)\cdot\left(g^0(x,y)+g^1(x)\right) \,dy\,dx,
$$
where $\omega$ is the frequency. The solution to the Euler-Lagrange equation for the resulting functional 
is a coupled system of equations for $g^0,$ $g^1,$ so that when the equation for $g^0$ is solved in terms of $g^1$ and substituted into the second equation, it takes the form (away from the sources):
\begin{equation}
{\mathcal A}^{\rm hom}g^1=\omega^2\beta(\omega^2)g^1,
\label{spectral_new}
\end{equation}
for some non-negative self-adjoint differential operator ${\mathcal A}^{\rm hom}$ and a special nonlinear function $\beta,$ which takes positive and negative values on alternating intervals of the real axis (leading to ``lacunae'', or ``band gaps'' in the spectrum of the corresponding operator) and is obtained from the spectral decomposition of $g^0$ and the subsequent averaging over $Y_0,$ see \cite{Zhikov}. From this point of view, the non-quadratic finite-strain functional $I_{\rm finite}$ is a ``matching'', "partially quadratic", homogenised model corresponding, {\it e.g.}, to finite-amplitude, rather than small-amplitude, wave motions that can no longer be treated using a quadratic model such as $I_{\rm small}$ but can still be used in place of  models of nonlinear elasticity where the elastic energy terms on both components of the composite (stiff and soft) are non-quadratic. 

\paragraph{Methods and previous results.}
In this paper we appeal to analytic methods that have been 
developed in the last two decades in the areas of nonlinear elasticity
and homogenisation. 
 Among these are the notion of two-scale convergence
introduced in \cite{Nguetseng}, \cite{Allaire}) and periodic unfolding 
(see \cite{Cioranescu_Damlamian} and references therein). The convergence statements
of our main results are expressed in the language of $\Gamma$-convergence (see \cite{DalMaso} and references therein). In order to
treat the geometric nonlinearity of the considered functional, we make use of the geometric rigidity
estimate (see \cite{FJM2002}). Since we consider a low energy regime,
linearisation and homogenisation take
place at the same time. The simultaneous treatment of both effects is inspired by
recent works \cite{Muller-Neukamm-10}, \cite{NeukammPhD2010}, \cite{Neukamm-12}, \cite{NeuVel}, 
\cite{HNV} of the third
author, where various problems involving \textit{simultaneous homogenisation, linearisation and dimension reduction} are studied. The
homogenisation of the kind of high-contrast composites that we study is related to the homogenisation for periodically perforated domains (e.g. see \cite{Oleinik1992},
\cite{Braides}). For instance, we make use of extensions across the
pores. As a side result we prove a version of the geometric
rigidity estimate for perforated domains (see Lemma~\ref{AL2} below). We would like to remark that while the present work is one of the few papers, along with \cite{CherCher}, \cite{Braides_ChiadoPiat_Piatnitski}, that treat the fully nonlinear high-contrast case, during the last decade there has been a significant amount of literature devoted to the mathematical analysis of phenomena associated with, or modelled by, a high degree of contrast between the properties 
of the materials constituting a composite, in the linearised setting. The first contributions in this direction are due to Zhikov \cite{Zhikov}, and Bouchitt\'{e} and Felbacq 
\cite{Bouchitte_Felbacq}, following an earlier paper by Allaire \cite{Allaire} and the collection of papers by Hornung {\it et al.} \cite{Hornung} (see 
also the references therein), where the special role of high-contrast
elliptic PDE was pointed out albeit not studied in detail.   These
works demonstrated that the behaviour of the field variable in such
models is of a two-scale type in the homogenisation limit, {\it i.e.} the
limit model cannot be reduced to a one-scale formulation and fields that
depend on the fast variable remain in the effective model.  They also noticed that the spectrum of such materials has a band-gap 
structure, as in (\ref{spectral_new}), and indicated how this fact could be exploited for high-resolution imaging and cloaking. 
It has since been an adopted approach to the theoretical construction of  ``negative refraction'' media, or more generally ``metamaterials'', which is now a hugely popular area of research in 
physics (see {\it e.g.} \cite{Pendry} and references therein). On the analytical side, a number of further works followed, in particular \cite{Bellieud}, \cite{BellieudGruais}, \cite{Briane1}, \cite{Briane2}, \cite{Camar}, 
\cite{ChSmZh}, \cite{KamotskySm}, \cite{Cherdantsev}, \cite{Smyshlyaev}, \cite{Braides_Briane}, \cite{Braides_ChiadoPiat_Piatnitski_2015}, where various consequences of high contrast (or, mathematically speaking, the property of non-uniform ellipticity) in the underlying equations have been explored. Among these are the ``non-locality'' and ``micro-torsion'' 
effects in materials with high-contrast inclusions in the shape of fibres extending in one or more directions, the ``partial band-gap'' wave propagation due to the high degree of anisotropy of one of the constituent media, and the localisation of energy in high-contrast media with a defect (``photonic crystal fibres''),
 all of which can be thought of as examples of ``non-standard'', 
or ``non-classical'', behaviour in composites, which is not available
in the usual moderate-contrast materials. 
In the present paper we aim to develop further a rigorous high-contrast theory in the context
of finite elasticity, where the underlying model is nonlinear.

With this paper we continue the multiscale theme initiated in
\cite{CherCher}, where the regime of large deformation gradients in the soft component 
of the composite was considered. Let us emphasise two points that contrast our contribution to some earlier work within the related field. First, we note that, apart from \cite{CherCher}, \cite{Braides_ChiadoPiat_Piatnitski}, a number of other articles ({\it e.g.}  \cite{Bouchitte_Bellieud}, \cite{Briane_monotonic}, \cite{Braides_Briane}) have treated high-contrast periodic composites in the nonlinear context. However, the related results are of limited relevance to nonlinear elasticity, due to the convexity or monotonicity assumptions made in these works. In the present paper we study a class of functionals subject to the requirement of material fame indifference (see assumption (W1) in Section \ref{S:ass}), which makes our analysis fit the fully nonlinear elasticity framework, as opposed to the works mentioned. Second, as was discussed above, the analysis of composites with  ``soft'' inclusions within a ``stiff'' matrix cannot be reduced to a ``decoupled'' model where the perforated medium obtained by removing the inclusions is considered first and the displacement within the inclusions is found independently, which from the physics perspective can be viewed as a kind of resonance phenomenon; {\it cf.} (\ref{spectral_new}) in the linearisation regime, for which an inherent energy coupling, in the limit as $\varepsilon\to0,$ between the soft and stiff components of the composite is essential. On a related note, the proof of the key compactness statement (Lemma \ref{L1}) involves the simultaneous analysis of the displacements on the two components. We would also like to highlight the fact that in \cite{CherCher} the order of the relative scaling of the displacements on the soft and stiff components of the composite are assumed from the outset, while in the present work it is the result of the above compactness argument itself. 


\paragraph{organisation of the paper.} In Section~\ref{S:ass} we state the
assumptions on the geometry of the composite and the material law. In
Section~\ref{S:main} we present the main results, starting with
results regarding two-scale compactness, convergence results in the
small strain regime and finally the convergence result in the finite
strain regime. All proofs are presented in Section~\ref{S:proofs}.

\subsection{Notation} 

Here we list some notation that we use throughout the text. Additional items will be introduced whenever they are first used in the text.

\begin{itemize}
\item $d\geq 2$ is the (integer) dimension of the space occupied by the material.

\item $p\ge1$ is the exponent in the notation $L^p$ for a Lebesgue space.  

\item $Y:=(0,1)^d$ the reference period cell; 
$Y^0$ is an open Lipschitz set whose closure is contained in $Y,$ and  
$Y^1:=Y\setminus\overline{Y^0}$.

\item $\Omega$, $\Omega^0_\e$ and $\Omega^1_\e$ denote the reference
  domains of the composite, the set occupied by the pore material, and the domain occupied by
  the matrix material, respectively, see Section~\ref{S:ass} for precise definition.
  
\item Unless stated otherwise, all function spaces $L^2(\Omega)$,
  $H^1(\Omega)$, $H^1_0(\Omega)$, {\it etc.} consist of functions taking values in $\R^d.$
  
\item Function spaces whose notation contains subscript ``c'' consist of functions that vanish outside a compact set.
  
\item The function spaces $H^1_\#,$ $H^1_0(Y^0),$ and $\mathcal A(Y^0)$ are introduced in Section~\ref{compact}.
  
\item We write $\cdot$ and $:$ for the canonical inner products in
  $\R^d$ and $\R^{d\times d}$, respectively.

\item $\SO(d)$ denotes the set of rotations in $\R^{d\times d}$.

\item $\lesssim$ stands for $\leq$ up to a multiplicative constant that only
depends on $d,$ $Y^1,$ $\Omega,$ and on $p$ if applicable. 
  
\end{itemize}

\section{Geometric and constitutive setup}\label{S:ass}

\noindent {\bf The pore geometry.} 
 The set $Y^0$ defined above describes the ``pores'' contained within the cell $Y$. Note that $Y^1$ is an open, bounded, connected set with Lipschitz boundary. Therefore, to each $\varphi\in H^1(Y^1)$ we can associate (see {\it e.g.} \cite{Oleinik1992}) a unique
\textit{harmonic} extension $g^1$ characterised by
\begin{equation}\label{eq:ext}
g^1=\varphi\mbox{\  \ in\  }Y^1,\qquad\qquad \int_{Y^0}\nabla
  g^1:\nabla\zeta\,dy=0\ \ \ \ \ \forall\zeta\in H^1_0(Y^0).
\end{equation}
For this extension the inequality
\begin{equation}\label{eq:ext2}
 \bigl\|\nabla g^1\bigr\|_{L^2(Y^0)}\leq C\|\nabla
  \varphi\|_{L^2(Y^1)}
\end{equation}
holds with a constant $C$ that only depends on $Y^1$. 

For a given domain $\Omega\subset\R^d$ and $\varepsilon>0,$ we define 
the sets $\Omega^0_\e$ and $\Omega^1_\e$ as follows:
\begin{equation*}
  \Omega^0_\e:=\bigcup\left\{\,\e(\xi+Y^0)\,|\,\xi\in\Z^d,\,\e(\xi+Y)\subset\Omega\,\right\},\qquad \Omega^1_\e:=\Omega\setminus\overline{\Omega^0_\e}.
\end{equation*}
Note that by construction $\Omega^1_\e$ is a Lipschitz domain. In particular, it
is connected and $\partial\Omega\subset\partial\Omega^1_\e$. We denote by
$\chi_\varepsilon$ the indicator function of the set of pores:
\begin{equation*}
  \chi_\e(x):=\left\{
    \begin{aligned}
      &1,&&x\in\Omega^0_\e,\\
      &0,&&x\in\R^d\setminus\Omega^0_\e.
    \end{aligned}
\right.
\end{equation*}

\noindent{\bf The composite.} The two materials are described by (Borel measurable) energy
densities $W^i:\R^{d\times d}\to[0,+\infty]$, $i=0,1$. Unless stated
otherwise, we assume that for $i=0,1:$
\begin{itemize}
\item[(W1)] $W^i$ is frame-indifferent, {\it i.e.} $W^i(RF)=W^i(F)$ for all $R\in
  \SO(d)$ and all $F\in\R^{d{\times}d};$
\item[(W2)] The identity matrix $I\in\R^{d\times d}$ is a ``natural
  state'', {\it i.e.} $W^i(I)=0$, and $W^i$ is non-degenerate, {\it i.e.} for all   \begin{equation*}
    W^i(F) \geq c_0\,\dist^2\bigl(F,\SO(d)\bigr),\ \ \ \forall F\in\R^{d{\times}d},\ \ \ \ \ c_0>0.
  \end{equation*}
\item[(W3)] $W^i$ has a quadratic expansion at $I$, {\it i.e.} there
  exists a non-negative quadratic form $Q^i$ on $\R^{d\times d}$ and an increasing
  function $r^i:[0,\infty)\to[0,\infty]$ with $\lim_{s\downarrow
    0}r^i(s)=0$, such that
  \begin{equation}\label{WTaylor}
 \Bigl|W^i(I+G) - Q^i(G)\Bigr|\leq |G|^2 r^i\bigl(|G|\bigr)\qquad\forall G\in\R^{d\times d}.
  \end{equation}
 \end{itemize}
As shown in \cite[Lemma~2.7]{Neukamm-12} the quadratic form $Q^i$
associated with $W^i$ via (W3) satisfies 
\begin{equation}\label{eq:sym}
  c_1|\sym\,G|^2\leq Q^i(G)=Q^i(\sym\,G)\leq c_1^{-1}|\sym\,G|^2\qquad\forall G\in\R^{d\times d},\ \ \ c_1>0.
\end{equation}
In the finite strain regime we consider a different set of assumptions
for $W^0,$ which are listed in Section~\ref{S:finite}.

\noindent{\bf The scaling parameter $\gamma$.} Throughout the paper $\gamma>0$ denotes a fixed
scaling parameter. It is a quantitative measure of the relative contrast between the two components of the composite.
\smallskip

\noindent{\bf Energy functional.} We define the elastic energy as a functional of the displacement, as follows:
\begin{equation}\label{eq:energy}
  \mathcal E_\e(u)=\int_\Omega
  \Bigl(\e^{2\gamma}W^0(I+\nabla u)\chi_\e+W^1(I+\nabla u)(1-\chi_\e)\Bigr)\,dx,\qquad
   u\in H^1_0(\Omega,\R^d).
\end{equation}
\section{Main results}\label{S:main}

\subsection{Compactness and two-scale convergence}
\label{compact}

We first present an {\it a priori} estimate and a two-scale compactness
statement for sequences $\varphi_\e\in
H^1_0(\Omega)$ whose energy is equi-bounded in the sense that
\begin{equation}\label{eq:89}
  \limsup\limits_{\e\downarrow 0}\Phi_\e^\gamma(\varphi_\e)<\infty,
\end{equation}
where
\begin{equation*}
  \Phi_\e^\gamma(\varphi):=\int_\Omega \dist^2\bigl(I+\nabla \varphi(x),\SO(d)\bigr)\,\bigl(\e^{2\gamma}\chi_\e+(1-\chi_\e)\bigr)\,dx.
\end{equation*}
Note that, by virtue of the non-degeneracy assumption (W2) the functional $\Phi^\gamma_\e(\cdot)$ bounds below 
$\mathcal E_\e({\rm id}+\cdot),$ where ${\rm id}(x)=x,$ $x\in\Omega.$
As we shall see in the upcoming Lemma~\ref{L1}, the inequality \eqref{eq:89} implies that
the sequence $\varphi_\varepsilon$ is bounded in  $H^1(\Omega)$, and thus weakly converges (up
to extracting a subsequence) to a
limit displacement $\varphi\in H^1_{0}(\Omega)$. The same argument applies to the more general situation of a rescaled displacement $m_\e\varphi_\e$, in which case we naturally rescale the functional $\Phi^\gamma_\e$ as well, see Lemma~\ref{L2}. For our
purpose we require a precise understanding of the oscillations that
emerge along that limit. We achieve this by combining two concepts:
\begin{itemize}
\item We write a representation for $\varphi_\e$ in the
spirit of an asymptotic decomposition as $\e\downarrow0.$
\item We study the convergence properties of the terms in this decomposition by appealing to two-scale convergence.
\end{itemize}
\medskip

In the following lemma we address the first item above.
\begin{lemma}
  \label{L1}
  Let  $\varphi_\e\in H^1_{0}(\Omega)$ and $0<\e\leq 1$. 
  \begin{enumerate}[(a)]
  \item There exists a unique pair of functions $g^0_\e\in H^1_0(\Omega^0_\e)$
    and $g^1_\e\in H^1_{0}(\Omega)$ such that
      \begin{equation}
        \label{L1:eq:1}
        \begin{aligned}
          &(i)\ \ \ \,\varphi_\e =g^1_\e+\e^{1-\gamma} g^0_\e,\\
          &(ii)\ \ \int_{\Omega^0_\e}\nabla g^1_\e:\nabla\zeta=0\ \ \forall\zeta\in H^1_0(\Omega^0_\e).
        \end{aligned}
      \end{equation}
  \item There exists a positive constant $C$ that only depends on
    $\Omega,Y^0$ such that
    \begin{equation}
      \label{L1:eq:2}
      \bigl\|g^0_\e\bigr\|_{L^2(\Omega)}^2 + \bigl\|\e\nabla
        g^0_\e\bigr\|^2_{L^2(\Omega)}+\bigl\|g^1_\e\bigr\|^2_{H^1(\Omega)}\,\leq C\,\Phi_\e^\gamma(\varphi_\e),
    \end{equation}
    where $g^0_\e,$ $g^1_\e$ and $\varphi_\e$ are related to each other as in (a).
  \end{enumerate}
\end{lemma}
Next, we recall the definition of two-scale convergence from
\cite{Nguetseng} and \cite{Allaire}:

\begin{definition}
  \label{def:4}
  We say that a sequence $f_\e\in L^2(\Omega)$ weakly
  two-scale converges to $f\in L^2(\Omega\times Y)$ if the sequence $f_\e$
  is bounded and 
  \begin{equation*}
    \lim\limits_{\e\downarrow 0}\int_\Omega
    f_\e(x)\varphi(x, x/\e)\,dx=\iint_{\Omega\times Y}f(x,y)\varphi(x,y)\,dxdy
  \end{equation*}
  for all $\varphi\in C^\infty_\mathrm{c}\bigl(\Omega,C^\infty_\#(Y)\bigr),$ where $C^\infty_\#(Y)$ is the set of infinitely smooth $Y$-periodic functions on ${\mathbb R}^d.$ We say that a sequence $f_\e\in L^2(\Omega)$, strongly
  two-scale converges to $f\in L^2(\Omega\times Y)$, if the sequence $f_\e$ weakly
  two-scale converges to $f$ and one has $\|f_\e\|_{L^2(\Omega)}\to
  \|f\|_{L^2(\Omega\times Y)}$ as $\e\downarrow0.$ For vector-valued functions two-scale
  convergence is defined component-wise.
\end{definition}
For convenience we introduce the following shorthands:
\begin{eqnarray*}
  f_\e\to f_0\qquad&:\Leftrightarrow&\qquad \text{$f_\e$ strongly
    converges to $f_0$ in $L^2(\Omega)$},\\    
  f_\e\wto f_0\qquad&:\Leftrightarrow&\qquad \text{$f_\e$ weakly  converges to $f_0$ in $L^2(\Omega)$},\\    
  f_\e\wtto f\qquad&:\Leftrightarrow&\qquad \text{$f_\e$ weakly
    two-scale converges to $f$ in $L^2(\Omega\times Y)$},\\
  f_\e\stto f\qquad&:\Leftrightarrow&\qquad \text{$f_\e$ strongly
    two-scale converges to $f$ in $L^2(\Omega\times Y)$}.
\end{eqnarray*}
The upcoming lemma states a two-scale compactness result for the displacements
$g^0_\e$ and $g^1_\e$ that appear in the representation \eqref{L1:eq:1}. Due to the
differential constraint satisfied by $g^0_\e$, the corresponding two-scale limits automatically satisfy certain structural properties, which can be
captured with the help of the following function spaces:
\begin{itemize}
\item $H^1_\#$ is the space of $[0,1)^d$-periodic
  functions in $H^1_\mathrm{loc}(\R^d)$.
\item $H^1_0(Y^0)$ is the closed
  subspace of  $H^1_\#$ consisting of functions $\psi\in H^1_\#$ with $\psi=0$ on $Y^1$.
\item $\mathcal A(Y^0)$ is the closed
  subspace of $H^1_\#$ consisting of functions $\psi\in H^1_\#$ that satisfy the identity
  \begin{equation*}
    \int_{Y^0}\nabla_y\psi:\nabla_y\zeta\,dy=0\qquad\forall\zeta\in H^1_0(Y^0).
  \end{equation*}
\end{itemize}
\begin{lemma}
  \label{L2}
  Consider a sequence $\varphi_\e\in H^1_0(\Omega)$ and let $(g^0_\e,g^1_\e)$ be associated 
  with $\varphi_\e$ via
  \eqref{L1:eq:1}. Suppose that there exists a sequence of positive
  numbers $m_\e$ such  that 
  \begin{equation*}
    \limsup\limits_{\e\downarrow 0}m_\e^{-2}\Phi_\e^\gamma(m_\e\varphi_\e)<\infty\qquad\mbox{
      and }\qquad m_\e=O(\e^\gamma).
  \end{equation*}
  Then there exist
  \begin{equation}
    g^0\in L^2\bigl(\Omega,H^1_0(Y^0)\bigr),\qquad g^1\in H^1_0(\Omega),\qquad
    \psi\in L^2\bigl(\Omega,\mathcal A(Y^0)\bigr)
  \end{equation}
  such that, up to selecting a subsequence, one has
  \begin{equation}\label{L2:2}
    \begin{split}
      &g^0_\e\wtto g^0,\qquad \e\nabla g^0_\e\wtto\nabla_y g^0,\\
      &g^1_\e\wto g^1\mbox{ weakly in $H^1(\Omega)$ and }\nabla g^1_\e\wtto \nabla
      g^1+\nabla_y\psi.
    \end{split}
  \end{equation}
\end{lemma}
The identification obtained in the previous lemma is sharp, in the sense of the following statement.
%
%
\begin{lemma}\label{L4}
  Let $g^0\in L^2\bigl(\Omega,H^1_0(Y^0)\bigr)$, $g^1\in H^1_0(\Omega)$ and
  $\psi\in L^2\bigl(\Omega,\mathcal A(Y^0)\bigr)$. Let $c_\e$ be an arbitrary sequence of positive numbers converging to zero. Then there exist function sequences
  $g^0_\e\in H^1_0(\Omega^0_\e)$, $g^1_\e\in H^1_0(\Omega)$ such that
  $(g^0_\e,g^1_\e)$ is related to $\varphi_\e:=g^1_\e+\e^{1-\gamma}
  g^0_\e$ as in \eqref{L1:eq:1}, and
  \begin{equation}\label{L4:1}
    \begin{split}
      &g^0_\e\stto g^0,\qquad \e\nabla g^0_\e\stto\nabla_y g^0,\\
      &g^1_\e\wto g^1\mbox{ weakly in $H^1(\Omega)$ and }\nabla g^1_\e\stto\nabla g^1+\nabla_y\psi,\\[0.4em]
      &\limsup\limits_{\e\downarrow 0}c_\e\bigl\|\nabla\varphi_\e\bigr\|_{L^\infty(\Omega)}=0.
    \end{split}
  \end{equation}
\end{lemma}
Our main result is formulated in terms of 
the notion of convergence
described in the above lemmas. For convenience we use the
following notation:
\def\decdef{\,\stackrel{\eqref{L1:eq:1}}{:=}\,}
\begin{itemize}
\item Given $\varphi_\e\in H^1(\Omega)$ we write $g^1_\e+\e^{1-\gamma}
  g^0_\e\decdef\varphi_\e$, if $g^1_\e\in H^1(\Omega)$, $g^0_\e\in
  H^1_0(\Omega_\e^0)$, and both functions are related to $\varphi_\e$ as in \eqref{L1:eq:1}.
\item We write $\varphi_\e\wtto (g^0,g^1)$, if $g^1_\e+\e^{1-\gamma}g^0_\e\decdef\varphi_\e$ and 
  \begin{equation}\label{eq:conv}
    g^0_\e\wtto g^0,\qquad\e\nabla g^0_\e\wtto\nabla_y g^0,\qquad g^1_\e\wto g^1\text{ weakly in }H^1(\Omega).
  \end{equation}
\item We write $\varphi_\e\stto (g^0,g^1)$, if  $g^1_\e+\e^{1-\gamma}
  g^0_\e\decdef\varphi_\e$ and 
  \begin{equation}\label{eq:convstrong}
    g^0_\e\stto g^0,\qquad \e\nabla g^0_\e\stto\nabla_y g^0,\qquad g^1_\e\wto g^1\text{ weakly in }H^1(\Omega).
  \end{equation}
\end{itemize}
\subsection{Convergence in the small strain regime
  $m_\e=o(\e^\gamma)$}\label{S:small}
\def\esmall{\textrm{small}}
\def\efin{\textrm{finite}}
Throughout this section we assume that the densities $W^0$ and $W^1$ satisfy the conditions (W1)--(W3).
We show that in the small strain regime the limit functional
\begin{equation*}
  \mathcal E_{\esmall}\,:\,L^2\bigl(\Omega,H^1_0(Y^0)\bigr)\times
  H^1_0(\Omega)\to[0,\infty),
\end{equation*}
is given by
\begin{equation*}
  \mathcal E_{\esmall}(g^0,g^1):=\int_{\Omega\times
    Y}\Bigl(Q^0\bigl(\nabla_yg^0(x,y)\bigr)+Q^1_{\ho}\bigl(\nabla g^1(x)\bigr)\Bigr)\,dx
\end{equation*}
where
\begin{equation}\label{def:Q1hom}
  Q^1_{\ho}(F):=\min\limits_{\psi\in \mathcal A(Y^0)}\int_{Y^1}Q^1\bigl(F+\nabla_y\psi(y)\bigr)\,dy.
\end{equation}
 More precisely, the following theorem holds.
\begin{theorem}
  \label{T1} Let $m_\e$ be a sequence of positive numbers and assume that
  $m_\e=o(\e^\gamma)$ as $\e\downarrow0.$
  \begin{enumerate}[(a)]
  \item (Compactness). Suppose that $\varphi_\e\in H^1_0(\Omega)$ satisfy
    \begin{equation*}
      \limsup\limits_{\e\downarrow 0}m_\e^{-2}\mathcal
      E_\e(m_\e\varphi_\e)<\infty.
    \end{equation*}
    Then, up to a subsequence, one has    $\varphi_\e\wtto (g^0,g^1)$ for some $g^0\in L^2\bigl(\Omega,H^1_0(Y^0)\bigr)$ and $g^1\in
    H^1_0(\Omega)$.

  \item (Lower bound). Consider $\varphi_\e\in H^1_0(\Omega)$ and  suppose
    that $\varphi_\e\wtto (g^0,g^1)$ for some $g^0\in L^2\bigl(\Omega,H^1_0(Y^0)\bigr)$ and $g^1\in
    H^1_0(\Omega)$. Then the estimate
    \begin{equation*}
      \liminf\limits_{\e\downarrow 0}m_\e^{-2}\mathcal
      E_\e(m_\e\varphi_\e)\geq \mathcal E_\mathrm{small}(g^0,g^1)
          \end{equation*}
  holds.
  \item (Recovery sequence). For all $g^0\in L^2\bigl(\Omega,H^1_0(Y^0)\bigr)$ and $g^1\in
    H^1_0(\Omega)$ there exists a sequence $\varphi_\e\in
    H^1_0(\Omega)$ such that $\varphi_\e\stto (g^0,g^1)$ and
    \begin{equation*}
      \lim\limits_{\e\downarrow 0}m_\e^{-2}\mathcal
      E_\e(m_\e\varphi_\e)=\mathcal E_\mathrm{small}(g^0,g^1).
    \end{equation*}
  \end{enumerate}
\end{theorem}

In the next result we consider a minimisation problem that involves
the density of the ``body forces" $\ell_\e\in L^2(\Omega)$. We study the variational limit of the
(scaled) total energy
\begin{equation}\label{eq:Ivarepsilon}
   I_\e(\varphi):=\frac{1}{m_\e^2}\left(\mathcal E_\e(m_\e\varphi)-\int_\Omega\ell_\e\cdot(m_\e\varphi)\,dx\right)
\end{equation}
where the scaling factor $m_\e$ is determined by the body forces via
\begin{subequations}
  \begin{equation}\label{eq:79}
    m_\e:=\e^{1-\gamma}\|\ell_\e\|_{L^2(\Omega^0_\e)}+\|\ell_\e\|_{L^2(\Omega)}.
  \end{equation}
In the small strain regime we assume that the body forces are small in
  the sense that
  \begin{equation}\label{eq:ass-smallstrain}
    m_\e=o(\e^\gamma), \ \ \ \e\downarrow0.
  \end{equation}
  Moreover, we assume that the (scaled) body-force densities converge, as $\e\downarrow0,$ in the following way:
  \begin{equation}\label{eq:90}
    m_\e^{-1}\e^{1-\gamma}\chi_\e\ell_\e\stto \ell^0,\qquad
    m_\e^{-1}\ell_\e\wto \ell^1.
  \end{equation}
\end{subequations}
It follows from Theorem \ref{T1}  that the variational limit of the total energy (\ref{eq:Ivarepsilon}) is given by the functional
\begin{equation}
\label{Ismall}
  I_{\esmall}(g^0,g^1):=\mathcal E_{\esmall}(g^0,g^1)-\int_\Omega\left(\int_{Y^0}\ell^0\cdot
      g^0\,dy+\ell^1\cdot g^1\right)\,dx.
\end{equation}
Next, we address convergence of infima and almost minimizers.
\begin{proposition}
  \label{P1}
  Assume that \eqref{eq:79}--\eqref{eq:90} hold.
  \begin{enumerate}[(b)]
  \item (Convergence of infima). One has
    \begin{equation*}
      \lim\limits_{\e\downarrow 0}\inf\limits_{\varphi\in
        H^1_0(\Omega)}I_\e(\varphi)=\min I_\mathrm{small}(g^0,g^1),
    \end{equation*}
    where the minimum on the right-hand side is taken over all
    $g^0\in L^2\bigl(\Omega,H^1_0(Y^0)\bigr)$ and  $g^1\in H^1_0(\Omega)$.
    Moreover, the minimum is attained for a unique pair $(g^0_*,g^1_*)$.
  \item (Convergence of minimisers). Let $\varphi_\e\in H^1_0(\Omega)$ be a sequence of almost
    minimisers, i.e.~
    \begin{equation}\label{P1:ass2}
      I_\e(\varphi_\e)\leq \inf\limits_{\varphi\in
        H^1_0(\Omega)}I_\e(\varphi) + o(1),\ \ \e\downarrow0.
    \end{equation}
    Then 
    \begin{equation*}
      \varphi_\e\wtto (g^0_*,g^1_*)\qquad\mbox{and}\qquad \nabla g^1_\e\wtto \nabla g^1_*+\nabla_y\psi_*
    \end{equation*}
    where $\psi_*\in L^2\bigl(\Omega,\mathcal A(Y^0)\bigr)$ denotes the unique ``corrector'' characterised by
    \begin{equation}\label{eq:91}
      Q_\ho^1\bigl(\nabla g^1_*(x)\bigr)=\int_{Y^1}Q^1\bigl(\nabla
      g^1_*(x)+\nabla_y\psi_*(x,y)\bigr)\,dy,\ \ \ \ \ \ \ \ \int_{Y}\psi_*(x,y)\,dy=0,
    \end{equation}
    for almost every $x\in\Omega$.
  \end{enumerate}
\end{proposition}
Next, we prove that almost minimisers
$\varphi_\e$ satisfy the asymptotic relation
\begin{equation}\label{eq:smallstrain-expansion}
  \varphi_\e=g^1_{*,\e}(x)+\e^{1-\gamma}g^0_{*,\e}(x)+o(1),\ \e\downarrow0,\qquad\text{in }W^{1,p}(\Omega)\qquad(p<2),
\end{equation}
where $g^1_{*,\e}$ and $g^0_{*,\e}$ formally obey the ``ansatz''
\begin{equation}\label{eq:93}
  g^0_{*,\e}(x)\stackrel{\text{formally}}{=} g^0_*(x,x/\e),\qquad g^1_{*,\e}(x)\stackrel{\text{formally}}{=}
  g^1_*(x)+\e\psi_*(x,x/\e),
\end{equation}
and the exponent $p$ depends on the regularity of the microstructure, see Assumption~\ref{ass:reg} below.
Here $(g^0_*,g^1_*)$ and $\psi_*$ denote the minimising pair and corrector
from Proposition~\ref{P1}. Since the functions on the right-hand
sides in \eqref{eq:93} are in general not smooth enough to define
$g^0_{*,\e}$ and $g^1_{*,\e}$ by \eqref{eq:93} directly, we use
instead the approximation associated with $(g^0_*,g^1_*,\psi_*)$ via Lemma~\ref{L4}.
\medskip

In addition to the properties of $Y^0$ assumed in
Section~\ref{S:ass}, we require the following assumption on the regularity of $Y^0$:
\begin{assumption}\label{ass:reg}
  There exist an exponent $p<2$ and a constant $C$ such that for all
  $\varphi\in H^1(Y^1)$ and $g^1\in H^1(Y)$ related via
  \eqref{eq:ext} we have $\|\nabla g^1\|_{L^p(Y^0)}\leq C\|\nabla
    \varphi\|_{L^p(Y^1)}$.
\end{assumption}
Note that Assumption~\ref{ass:reg} is satisfied if $Y^0$ can be written as
the disjoint union of a finite number of Lipschitz domains
$Y^0_1,\ldots Y^0_N$ with $\partial Y^0_i\cap\partial
Y^0_j=\emptyset$ for $i\neq j$.
\begin{theorem}\label{T1b}
  Assume that \eqref{eq:79}--\eqref{eq:90} hold, and let Assumption~\ref{ass:reg} be satisfied. Let $\varphi_\e\in H^1_0(\Omega)$ be a sequence of almost minimisers, i.e.~
  \begin{equation*}
    I_\e(\varphi_\e)\leq \inf\limits_{\varphi\in H^1_0(\Omega)}I_\e(\varphi) + o(1),\ \ \e\downarrow0.
  \end{equation*}
  Let $(g^0_*,g^1_*)$ be the minimiser of $I_\mathrm{small}$ and let
  $\psi_*$ be defined through \eqref{eq:91}. Let
  $(g^0_{*,\e},g^1_{*,\e})$ and
  $\varphi_{*,\e}=g^1_{*,\e}+\e^{1-\gamma} g^0_{*,\e}$ be associated
  with $(g^0_*,g^1_*,\psi_*)$ as in  Lemma~\ref{L4}, {\it i.e.} $\varphi_{*,\e}\stto (g^0_*,g^1_*)$ and $\nabla g_{*,\e}^1\stto \nabla g^1_*+\nabla_y\psi_*.$  Then for $g^1_\e+\e^{1-\gamma}g^0_\e\decdef\varphi_\e$ one has
  \begin{equation}\label{T1b:1}
   \bigl\|g^0_\e-g^0_{*,\e}\bigr\|_{L^p(\Omega_\e^0)}
    +\bigl\|\e\nabla g^0_\e-\e\nabla g^0_{*,\e}\bigr\|_{L^p(\Omega_\e^0)}
    +\bigl\|g^1_\e- g^1_{*,\e}\bigr\|_{W^{1,p}(\Omega)}\to 0\ \ \ \mathrm{as\ }\e\downarrow0.
  \end{equation}
\end{theorem}

\begin{remark}
  \label{sec:conv-small-stra}
  To illustrate the result of Theorem~\ref{T1b}, consider the case
  $\gamma=1$ with $\ell_\e:=m_\e\ell(x,\tfrac{x}{\e}),$ where
  $\ell(x,y)$ is smooth both in $x$ and $y$ and is periodic in $y$. If the domain $\Omega$ and the pore
  set $Y^0$ are sufficiently regular, 
the minimisers
  $(g^0_*,g^1_*)$ and $\psi_*$ are smooth, by the classical elliptic regularity theory, 
see {\it e.g.} \cite{Evans}. In that case we may set
  \begin{equation*}
    g^0_{*,\e}(x):=g^0_*(x, x/\e)\ \mbox{ and }\ g^1_{*,\e}:=g^1_*(x)+\e\psi_*(x, x/\e),
  \end{equation*}
  and the asymptotic formula for $\varphi_\varepsilon$ reads
  \begin{equation*}
    \varphi_\e=g^0_*(x, x/\e)+g^1_*(x) + \e\psi_*(x, x/\varepsilon)+R_\varepsilon(x),
  \end{equation*}
where $\Vert R_\varepsilon\Vert_{W^{1,p}(\Omega)}\to0$ as $\e\downarrow0.$  
\end{remark}

\begin{remark}[Example in Section \ref{intro_section}]\label{R:intro:small}
  Suppose that $\alpha>1$. If we specialize Theorem~\ref{T1},  Proposition~\ref{P1} and Theorem~\ref{T1b} to $\gamma=1$, $m_\e=\e^\alpha$ and $\ell_\e=\e^\alpha f$ for some $f\in L^2(\Omega)$, then we recover the special case (in the small strain regime) presented in the Introduction. In particular, Proposition~\ref{P1} proves the claimed two-scale $\Gamma$-convergence and Theorem~\ref{T1b} establishes a two-scale expansion. Next we argue that the functionals $I^\alpha_\e$ $\Gamma$-converge (with respect to the weak topology in $L^2(\Omega)$) to $\bar I_{\text{\rm small}}$ in \eqref{intro:small:weak-single-scale}. To this end, note that with the definition of $\overline Q$ from the introduction we get for arbitrary $(g^0,g^1)$ (with $\overline G:=\int_{Y_0}g^1\,dy)$ the inequality:
  \begin{eqnarray*}
    I_{\text{\rm small}}(g^0,g^1)&\geq&\int_{\Omega}Q^1_{\hom}\bigl(\nabla g^1(x)\bigr)-\bigl(g^1(x)+\overline G(x)\bigr)\cdot f(x)\,dx \\
                                 &&+\int_{\Omega}\min\left\{\int_{Y^0}Q^0\bigl(\nabla_y g^0(y)\bigr)\,:\,g^0\in H^1_0(Y^0),\,\int_{Y^0}g^0(y)=\overline G(x)\,\right\}\\[0.4em]
    &\geq & \int_{\Omega}\overline Q(\overline G(x))+Q^1_{\hom}\bigl(\nabla g^1(x)\bigr)-\bigl(g^1(x)+\overline G_*(x)\bigr)\cdot f(x)\,dx\\[0.4em]
                                 &\geq & \bar I_{\text{\rm small}}(g^1+\overline G).
  \end{eqnarray*}
  Together with part (b) of Theorem~\ref{T1} this proves the lower-bound part of the $\Gamma$-convergence statement. On the other hand, for arbitrary $\varphi\in L^2(\Omega)$ we can find a unique pair $(g^0,g^1)$ such that $\varphi=g^1+\int_{Y_0}g^0\,dy$ and $I_{\text{\rm small}}(g^0,g^1)=\bar I_{\text{\rm small}}(\varphi)$. Together with part (c) of Theorem~\ref{T1} we get existence of a recovery sequence. Hence, we conclude that $I^\alpha_\e$ $\Gamma$-converge (weakly in $L^2(\Omega)$) to $\bar I_{\text{\rm small}}$.

  In particular, since  $(g^0_*,g^1_*)$ minimises $I_{\textrm{small}}$, we get, as claimed, that
  \begin{eqnarray*}
    I_{\text{\rm small}}(g^0_*,g^1_*)&=& \bar I_{\text{\rm small}}(g_*^1+\overline G_*)=\bar I_{\text{\rm small}}(\varphi_*),
  \end{eqnarray*}
  where $\varphi_*=g^1_*+\overline G_*$ is the weak limit in $L^2(\Omega)$ of any sequence of almost minimisers of $I^\alpha_\e$.
\end{remark}

\subsection{Convergence in the finite strain regime
  $m_\e=\e^\gamma$}\label{S:finite}
Throughout this section we assume that 
\begin{itemize}
\item $W^1$ satisfies the conditions (W1)--(W3).
\item $W^0:\R^{d\times d}\to[0,\infty)$ is
continuous and satisfies the growth condition
\begin{subequations}
  \begin{equation}\label{eq:assW0a}
    c_0\,\dist^2\bigl(F,\SO(d)\bigr)\leq W^0(F)\leq c_0^{-1}\bigl(1+|F|^2\bigr)\qquad\forall F\in{\mathbb R}^{d\times d},
  \end{equation}
  and the local Lipschitz condition
  \begin{equation}\label{eq:assW0b}
   \Bigl|W^0(F+G)-W^0(F)\Bigr|\leq c_0^{-1}\bigl(1+|F|+|G|\bigr)|G|\qquad\forall F,G\in{\mathbb R}^{d\times d}.
  \end{equation}
\end{subequations} 
\end{itemize}
We prove that in the finite strain regime the limit functional
\begin{equation*}
  \mathcal E_{\efin}\,:\,L^2\bigl(\Omega,H^1_0(Y^0)\bigr)\times
  H^1_0(\Omega)\to[0,\infty)
\end{equation*}
is given by
\begin{equation*}
  \mathcal E_{\efin}(g^0,g^1):=\iint_{\Omega\times Y^0}\mathcal QW^0\bigl(I+\nabla_yg^0(x,y)\bigr)\,dydx+\int_\Omega Q^1_{\ho}\bigl(\nabla g^1(x)\bigr)\,dx,
\end{equation*}
where $\mathcal QW^0$ denotes the quasiconvex envelope of $W^0$ (see {\it e.g.} \cite{Dacorogna}). The associated limit of the total energy $I^\varepsilon,$ see \eqref{eq:Ivarepsilon}, is given by ({\it cf.} (\ref{Ismall}))
\begin{equation*}
  I_\mathrm{finite}(g^0,g^1):=\mathcal E_\mathrm{finite}(g^0,g^1)-\int_\Omega\left(\int_{Y^0}\ell^0\cdot
      g^0\,dy+\ell^1\cdot g^1\right)\,dx,
\end{equation*}
where $\ell^0,$ $\ell^1$ are defined in the same way as in (\ref{eq:90}).

\begin{theorem}
  \label{T2}
 \begin{enumerate}[(a)]
  \item (Lower bound). Consider a sequence $\varphi_\e\in
    H^1_0(\Omega)$ and the associated decomposition $\e^{1-\gamma}g^0_\e+g^1_\e\decdef\varphi_\e$. If $\varphi_\e\wtto (g^0,g^1)$ and 
    $g^0_\e\stto g^0$, then
    \begin{equation*}
      \liminf\limits_{\e\downarrow 0}\e^{-2\gamma}\mathcal
      E_\e(\e^\gamma\varphi_\e)\geq \mathcal E_\mathrm{finite}(g^0,g^1).
    \end{equation*}
  \item (Recovery sequence). For any $g^0\in L^2\bigl(\Omega,H^1_0(Y^0)\bigr)$ and $g^1\in
    H^1_0(\Omega)$ there exists a sequence $\varphi_\e\in
    H^1_0(\Omega)$ such that $\varphi_\e\wtto (g^0,g^1)$ and
    \begin{equation*}
      \lim\limits_{\e\downarrow 0}\e^{-2\gamma}\mathcal
      E_\e(\e^\gamma\varphi_\e)=\mathcal E_\mathrm{finite}(g^0,g^1).
    \end{equation*}
\item Suppose that the force densities $\ell_\e\in
  L^2(\Omega)$ satisfy \eqref{eq:79} and \eqref{eq:90} with $m_\e=\e^{\gamma}$. Then the infima converge,
  i.e.
  \begin{equation}\label{T2:1}
    \lim\limits_{\e\downarrow 0}\inf\limits_{\varphi\in
      H^1_0(\Omega)}I_\e(\varphi)=\inf I_\mathrm{finite}(g^0,g^1),
  \end{equation}
  where the infimum on the right-hand side is taken over all
  functions $g^0\in L^2\bigl(\Omega,H^1_0(Y^0)\bigr)$ and  $g^1\in H^1_0(\Omega)$.
  Moreover, there exist a
  minimising pair $(g^0_*,g^1_*)$ and a recovery sequence $\varphi_\e\in
  H^1_0(\Omega)$ with $\varphi_\e\stto (g^0_*,g^1_*)$ such
  that
  \begin{equation}\label{T2:2}
    I_\e(\varphi_\e)\to I_\mathrm{finite}(g^0_*,g^1_*)=\min I_\mathrm{finite}(g^0,g^1)\ \ \ \mathrm{as}\ \e\downarrow0.
  \end{equation}
  \end{enumerate}
\end{theorem}
\begin{remark}[Example in Section \ref{intro_section}]\label{R:intro:finite}
  If we consider Theorem~\ref{T2} and Proposition~\ref{P1} in the case $\gamma=1$, $m_\e=\e$ and $\ell_\e=\e f$ for some $f\in L^2(\Omega)$, then we recover the special case (in the finite strain regime) presented in the introduction. In particular, we deduce that the functionals $I^\alpha_\e$ two-scale $\Gamma$-converge (in the sense of Theorem~\ref{T2}) to $I_{\text{\rm finite}}$. Arguing as in the small-strain regime, cf. Remark~\ref{R:intro:small}, we deduce that $I^\alpha_\e$ $\Gamma$-converges (with respect to the weak topology in $L^2(\Omega)$) to $\bar I_{\text{\rm finite}}$, cf. \eqref{intro:finite:weak-single-scale}. 
\end{remark}

%



\section{Proofs}\label{S:proofs}

We start by proving the auxiliary results discussed in Section \ref{compact}. 
Sections \ref{S:P2} and \ref{finitestrain} contain the proofs of the main statements in the small strain and finite strain cases, respectively.

\subsection{Proofs of Lemma~\ref{L1}, Lemma~\ref{L2}, and
  Lemma~\ref{L4}:  a priori
  estimate, compactness and approximation}
A key ingredient in the proof of Lemma~\ref{L1} is the geometric
rigidity estimate by Friesecke~{\it et al.} \cite{FJM2002}:
\begin{theorem}[Geometric rigidity estimate, see \cite{FJM2002}]\label{T:rigidity}
  Let $U$ be an open, bounded Lipschitz domain in $\R^d$, $d\geq 2$. There exists a constant $C(U)$ with the following property: for each 
  $v\in H^1(U)$ there is a rotation $R\in\SO(d) $ such that
  
\begin{equation*}
    \int_U\bigl|\nabla v(x)-R\bigr|^2\,dx\leq
    C(U)\int_U\dist^2\bigl(\nabla v(x),\SO(d)\bigr)\,dx.
  \end{equation*}
  Moreover, the constant $C(U)$ is invariant under uniform scaling of $U$.
\end{theorem}
In fact, we need the following modified version, which is adapted to perforated domains.
\begin{lemma}\label{AL2}
  There exists a constant $C>0$ that only depends on $\Omega$ and $Y^1$ such that
  for all $\e>0$ and $v\in H^1(\Omega)$ satisfying
  \begin{equation}\label{AL2:ass}
    \int_{\Omega^0_\e}\nabla v:\nabla\zeta\,dx=0\qquad\forall\zeta\in H^1_0(\Omega^0_\e),
  \end{equation}
 the estimates
  \begin{eqnarray}
    \label{AL2:1}
    \bigl\|\dist\bigl(\nabla v, \SO(d)\bigr)\bigr\|_{L^2(\Omega)}
    &\leq& C\, \bigl\|\dist\bigl(\nabla {v}, \SO(d)\bigr)\bigr\|_{L^2(\Omega_\e^1)},\\
    \label{AL2:2}
    \|\nabla v - R\|_{L^2(\Omega_\e^1)} &\leq& C\,\bigl\|\mathrm{dist}
    \bigl(\nabla v, \SO(d)\bigr)\bigr\|_{L^2(\Omega_\e^1)},
  \end{eqnarray}
 hold for some $R\in\SO(d),$ which may depend on $v.$  In addition, if $v(x)=x+c$ on $\partial\Omega$ for some constant
  $c$, then we may set $R=I$.
\end{lemma}
\begin{proof}[Proof of Lemma~\ref{AL2}]
\smallskip

\step 1 The proof of the inequality \eqref{AL2:1}.

Let $\widehat\Omega_\e:=\bigcup\left\{\,\e(\xi+Y)\,|\,\xi\in\Z^d,\,\e(\xi+Y)\subset\Omega\,\right\}$
denote the union of $\e$-cells that are completely contained in
$\Omega$. Since $\Omega\setminus\widehat\Omega_\e\subset\Omega^1_\e$,
it suffices to  prove
\eqref{AL2:1} for $\Omega$ replaced by
$\widehat\Omega_\e$,
respectively. In fact we shall prove the following stronger estimate:   for
all $\xi\in\Z^d$ with $\e(\xi+Y)\subset\Omega$ we have 
\begin{equation}\label{AL2:pf:1}
  \int_{\e(\xi+Y)}\dist^2\bigl(\nabla v, \SO(d)\bigr)\,dx
  \lesssim\, \int_{\e(\xi+Y^1)}\dist^2\bigl(\nabla {v}, \SO(d)\bigr)\,dx.
\end{equation}
For the argument fix an admissible $\xi\in\Z^d$. Application of
Theorem~\ref{T:rigidity} with $U=\e(\xi+Y^1)$ yields a
rotation $R\in\SO(d)$ such that
\begin{equation}\label{AL2:pf:2}
  \int_{\e(\xi+Y^1)}|\nabla v-R|^2\,dx
  \lesssim\, \int_{\e(\xi+Y^1)}\dist^2\bigl(\nabla{v}, \SO(d)\bigr)\,dx.
\end{equation}
Note that the multiplicative constant in the estimate above only
depends on $Y^1$, since $\e(\xi+Y^1)$ is a dilation and translation of
$Y^1$. On the other hand, since  $\e(\xi+Y^0)\ni x\mapsto (v(x)-Rx)$ is harmonic, we have (cf. \eqref{eq:ext2}):
\begin{equation*}
  \int_{\e(\xi+Y)}\dist^2\bigl(\nabla v,\SO(d)\bigr)\,dx\leq
  \int_{\e(\xi+Y)}|\nabla v-R|^2\,\lesssim\,\int_{\e(\xi+Y^1)}|\nabla v-R|^2.
\end{equation*}
Combined with \eqref{AL2:pf:2}, inequality \eqref{AL2:pf:1} follows.
\smallskip
  
\step 2 The proof of the rigidity estimate \eqref{AL2:2}.

  From \eqref{AL2:1} and
  Theorem~\ref{T:rigidity} (applied with $U=\Omega$) we deduce that for some $R\in \SO(d)$:
  \begin{equation}\label{AL2:pf:0}
    \|\nabla v - R\|_{L^2(\Omega)}\,\lesssim\,\bigl\|\dist\bigl(\nabla v, \SO(d)\bigr)\bigr\|_{L^2(\Omega^1_\e)},
  \end{equation}
  which in particular implies \eqref{AL2:2}. Finally we argue
  that one can set $R=I$, if $v=x+c$ on
  $\partial\Omega$. In view of \eqref{AL2:pf:0}, it suffices to show
  that $\int_\Omega|\nabla v-I|^2\,dx\leq\int_\Omega|\nabla
    v-R|^2\,dx$  for all $R\in\SO(d)$. This inequality can be seen as
    follows: Consider $\varphi(x):=v(x)-x-c$
  and note that $\varphi$ vanishes on $\partial\Omega$, so that 
  \begin{equation*}
    \int_\Omega|\nabla v-I|^2\,dx=\int_\Omega|\nabla\varphi|^2\,dx\leq
    \int_\Omega|\nabla\varphi|^2+|I-R|^2\,dx=\int_\Omega\bigl|\nabla\varphi+(I-R)\bigr|^2\,dx=\int_\Omega|\nabla
    v-R|^2\,dx,
  \end{equation*}
which in fact holds for an arbitrary matrix $R$.
\end{proof}

We are now in position to present the proofs of Lemmas \ref{L1} and \ref{L2}.

\begin{proof}[Proof of Lemma~\ref{L1}]
  In the following, the symbol  $\lesssim$ stands for $\leq$ up to a multiplicative
  constant that only depends on $Y^1$ and $\Omega$.
  \smallskip
  
  \step 1 Existence of the decomposition (\ref{L1:eq:1}) and derivation of 
the estimate for $g^1_\e$.
  
  Let $g^1_\e$ denote the unique function in $H^1(\Omega)$
  characterised by $g^1_\e=\varphi_\e$ in $\Omega^1_\e$ and
  \eqref{L1:eq:1} (ii). Since $\partial\Omega^0_\e$ is Lipschitz, we
  deduce that $g^0_\e:=\e^{\gamma-1}(\varphi_\e-g^1_\e)\in
  H^1_0(\Omega^0_\e)$. This proves the existence of the decomposition.
  We claim that
  \begin{equation}\label{eq:72}
    \int_{\Omega^0_\e}\bigl|\nabla
    g^1_\e\bigr|^2\,dx\,\lesssim\,\int_{\Omega^1_\e}|\nabla\varphi_\e|^2\,dx=\int_{\Omega^1_\e}\bigl|\nabla
    g^1_\e\bigr|^2\,dx.
  \end{equation}
  Since $\Omega^0_\e$ is defined as the union of the sets 
  $\e(\xi+Y^0)$ with $\xi\in
  Z_\e:=\{\,\xi\in\Z^d\,:\,\e(\xi+Y)\subset\Omega\,\}$, it
  suffices to prove $\int_{\e(\xi+Y^0)}\bigl|\nabla
  g^1_\e\bigr|^2\,dx\,\lesssim\,\int_{\e(\xi+Y^1)}|\nabla\varphi_\e|^2\,dx$.
  The latter follows from \eqref{eq:ext2} by a scaling argument, since
  the rescaled functions $y\mapsto\varphi_\e\bigl(\e(\xi+y)\bigr)$ and $y\mapsto
  g^1_\e\bigl(\e(\xi+y)\bigr)$ satisfy \eqref{eq:ext}.

  Next, we prove \eqref{L1:eq:2}. Consider
  $v_\e(x):=x+g^1_\e(x)$ and note that $v_\e(x)$ satisfies
  \eqref{AL2:ass}. Hence, \eqref{eq:72} and Lemma~\ref{AL2} yield
  \begin{equation}\label{L1:pf:1}
    \int_{\Omega}\bigl|\nabla g^1_\e\bigr|^2\,dx\lesssim \int_{\Omega_\e^1}\bigl|\nabla g^1_\e\bigr|^2\,dx\leq\int_\Omega|\nabla
    v_\e-I|^2\stackrel{\eqref{AL2:2}}{\leq}
    C\int_{\Omega^1_\e}\dist^2\bigl(I+\nabla\varphi_\e(x),\SO(d)\bigr)\,dx\leq \Phi^\gamma_\e(\varphi_\e).
  \end{equation}
  Since $g^1_\e$ vanishes on the boundary of $\Omega$, the estimate
  upgrades (by Poincar\'e's inequality) to $\bigl\|g^1_\e\bigr\|_{H^1(\Omega)}\leq
  C\Phi^\gamma_\e(\varphi_\e)$.
  \smallskip

  \step 2 Derivation of the estimate for $g^0_\e$.

  Since we have an improved Poincar\'e inequality (see {\it e.g.} \cite[Lemma~1.6]{Hornung}):
  \begin{equation}\label{improvedPoincare}
    \forall g\in H^1_0(\Omega^0_\e)\,:\qquad
    \|g\|_{L^2(\Omega^0_\e)}\lesssim\e\|\nabla g\|_{L^2(\Omega^0_\e)},
  \end{equation}
  it suffices to prove
  \begin{equation*}
  \bigl\|\e\nabla g^0_\e\bigr\|_{L^2(\Omega^0_\e)}^2\lesssim \Phi^\gamma_\e(\varphi_\e).
  \end{equation*}
  To this end, notice that since $\varphi_\e$
  vanishes on the boundary of $\Omega$, we have 
  \begin{equation}\label{L1:pf:3}
    \begin{aligned}
      \|\nabla\varphi_\e\|^2_{L^2(\Omega)}\,=\,&\min_{R\in\SO(d)}\|I+\nabla\varphi_\e-R\|^2_{L^2(\Omega)}\stackrel{\text{Theorem~\ref{T:rigidity}}}{\lesssim}
      \int_\Omega\dist^2\bigl(I+\nabla\varphi_\e(x),\SO(d)\bigr)\,dx\\
      \leq\,&\e^{-2\gamma}\Phi^\gamma_\e(\varphi_\e).
    \end{aligned}
  \end{equation}
  Thanks to the first identity in \eqref{L1:eq:1}, we get by triangle inequality:
  \begin{equation*}
    \bigl\|\e^{1-\gamma}\nabla g^0_\e\bigr\|_{L^2(\Omega^0_\e)}\leq\|\nabla\varphi_\e\|_{L^2(\Omega)}+\bigl\|\nabla g^1_\e\bigr\|_{L^2(\Omega)}.
  \end{equation*}
  Combined with \eqref{L1:pf:1} and \eqref{L1:pf:3} we finally get
  \begin{equation*}
    \bigl\|\e\nabla
    g^0_\e\bigr\|^2_{L^2(\Omega^0_\e)}=\e^{2\gamma}\bigl\|\e^{1-\gamma}\nabla
    g^0_\e\bigr\|^2_{L^2(\Omega^0_\e)}\lesssim\Phi^\gamma_\e(\varphi_\e).
  \end{equation*}
\end{proof}

\begin{proof}[Proof of Lemma~\ref{L2}]

  \step 1 {\it A priori} estimate and basic compactness.
  
  From  Lemma~\ref{L1} (applied with $\varphi_\e$ replaced by $m_\e\varphi_\e$) we deduce that
  \begin{equation}
    \label{L2:pf:1}
    \limsup\limits_{\e\downarrow 0}\left(\bigl\|g^0_\e\bigr\|_{L^2(\Omega)}^2 + \bigl\|\e\nabla
    g^0_\e\bigr\|^2_{L^2(\Omega)}+\bigl\|g^1_\e\bigr\|^2_{H^1(\Omega)}\right)\,\leq C
  \,\limsup\limits_{\e\downarrow 0}m_\e^{-2}\Phi_\e^\gamma(m_\e\varphi_\e)<\infty.
  \end{equation}
  Hence, by standard results concerning two-scale convergence ({\it cf.} \cite[Proposition 1.14]{Allaire} 
  and \cite[Proposition 4.2]{Visintin}), there exist $g^1\in
  H^1_0(\Omega)$, $\psi\in L^2\bigl(\Omega,H^1_\#\bigr)$ and 
  $g^0\in L^2\bigl(\Omega, H^1_\#\bigr)$ such that, up to a subsequence, one has 
  \begin{align*}
    &g^1_\e\wto g^1\mbox{ weakly in }H^1(\Omega), \qquad \nabla g^1_\e\wtto \nabla g^1+\nabla_y\psi,\\[0.4em]
    &g^0_\e\wtto g^0, \qquad \e\nabla g^0_\e\wtto \nabla_yg^0.
  \end{align*}
  
  \step 2 The proof of the inclusion  $\psi\in L^2\bigl(\Omega,\mathcal A(Y^0)\bigr)$.

  By a density argument, it suffices to show that
  \begin{equation}\label{L2:pf:2}
    \iint_{\Omega\times Y^0}\nabla_y\psi(x,y):\nabla_y\bigl(\zeta_1(x)\zeta_2(y)\bigr)\,dxdy=0
  \end{equation}
  for all scalar functions $\zeta_1\in C^\infty_\mathrm{c}(\Omega)$, and all $\zeta_2\in C^\infty_\mathrm{c}(Y^0)$.
  To this end, we identify $\zeta_2$ with its unique $Y$-periodic
  extension to $\R^d$ that vanishes on $Y^1$, and set
  \begin{equation*}
    \zeta_\e(x):=\e\zeta_1(x)\zeta_2(x/\e),\ \ \ x\in\Omega.
  \end{equation*}
  Thanks to \eqref{L1:eq:1} we have
  \begin{equation*}
        \int_\Omega\nabla g^1_\e:\nabla\zeta_\e\,dx=0.
  \end{equation*}
  As can be easily checked, we have  $\nabla\zeta_\e\stto \zeta_1(x)\nabla_y\zeta_2(y)$, so that
  \begin{eqnarray*}
    0&=&\lim\limits_{\e\downarrow 0}\int_\Omega\nabla
    g^1_\e:\nabla\zeta_\e\,dx\\
    &=&\iint_{\Omega\times Y^0}\bigl(\nabla g^1(x)+\nabla_y\psi(x,y)\bigr):\bigl(\zeta_1(x)\nabla_y\zeta_2(y)\bigr)\,dxdy\\
    &=&\iint_{\Omega\times Y^0}\nabla_y\psi(x,y):\bigl(\zeta_1(x)\nabla_y\zeta_2(y)\bigr)\,dxdy,
  \end{eqnarray*}
where the last identity holds thanks to the periodicity of $\zeta_2$. This proves \eqref{L2:pf:2}.
  \smallskip

  \step 3 The proof of the inclusion $g^0\in L^2\bigl(\Omega,H^1_0(Y^0)\bigr)$.

  By a density argument, it suffices to show that
  \begin{equation}\label{L2:pf:3}
    \iint_{\Omega\times Y}g^0(x,y)\cdot\bigl(\zeta_1(x)\zeta_2(y)\bigr)\,dxdy=0
  \end{equation}
  for all scalar functions $\zeta_1\in C^\infty_\mathrm{c}(\Omega)$ and all $\zeta_2\in H^1_\#$ with $\zeta_2=0$ on $Y^0$. We argue by considering the function $\zeta_\e(x):=\zeta_1(x)\zeta_2(\tfrac{x}{\e}),$ $x\in\Omega,$
  the support of which is contained in $\Omega^1_\e$ for $\e\ll 1$.
  Since $\zeta_\e\stto \zeta_1(x)\zeta_2(y)$, and since $g^0_\e$ is
  supported in $\Omega^0_\e$, we deduce that 
  \begin{equation*}
    0=\lim\limits_{\e\downarrow 0}\int_\Omega
    g^0_\e(x)\cdot\zeta_\e(x)\,dx=\iint_{\Omega\times Y} g^0(x,y)\cdot\bigl(\zeta_1(x)\zeta_2(y)\bigr)\,dxdy.
  \end{equation*}
  This completes the argument.
\end{proof}
In the proof of Lemma~\ref{L4} we appeal to the construction of a
diagonal sequence that is due to Attouch, see \cite{Attouch1984}:
\begin{lemma}
\label{lem:attouch} 
  For any
  $h:[0,\infty)^2\rightarrow
  [0,+\infty],$ there exists a mapping $(0,1)\ni\varepsilon\mapsto
  \delta(\varepsilon)\in(0,1)$ such that
  \begin{equation*}
    \lim\limits_{\varepsilon\downarrow
      0}\delta(\varepsilon)=
    0\qquad\text{and}\qquad\limsup\limits_{\varepsilon\downarrow
      0}h\bigl(\e,\delta(\varepsilon)\bigr)\leq
    \limsup\limits_{\delta\downarrow
      0}\limsup\limits_{\varepsilon\downarrow 0}h(\e,\delta).
  \end{equation*}
\end{lemma}

\begin{proof}[Proof of Lemma~\ref{L4}]
\step 1 Characterisation of strong two-scale convergence via unfolding.

For $f_\e:\Omega\to\R$ and $f:\Omega\times Y\to\R$ define
\begin{equation*}
  d_\e(f_\e,f):=\int_{\R^d}\int_Y\bigl|\tilde f_\e(\e\lfloor x/\e\rfloor+\e y)-\tilde f(x,y)\bigr|^2\,dydx
\end{equation*}
where $\tilde f_\e$ denotes the extension by zero of $f_\e$ to $\R^d$,
$\tilde f$ denotes the extension by zero of $f$ to $\R^d\times Y$, and
$\lfloor z\rfloor$ denotes the unique element in $\Z^d$ with $z-\lfloor
z\rfloor\in[0,1)^d$. We recall from \cite{Visintin} that
\begin{equation}
  \label{eq:96}
  f_\e\stto f\qquad\Leftrightarrow\qquad d_\e(f_\e,f)\to 0.
\end{equation}
The characterisation extends in the obvious way to vector-valued  functions.

\step 2 Construction of $g^1_\e$.

We claim that there exists a sequence $g^1_\e$ in $H^1_0(\Omega)$ whose elements satisfy
\eqref{L1:eq:1}(ii) and
\begin{equation}\label{L4:pf0}
  g^1_\e\wto g^1\mbox{ weakly in }H^1(\Omega),\qquad \nabla g^1_\e\stto g^1+\nabla_y\psi,\qquad
  \limsup\limits_{\e\downarrow 0}c_\e\bigl\|\nabla g^1_\e\bigr\|_{L^\infty(\Omega)}=0.
\end{equation}
Indeed, by a density argument there exist $g^{1,\delta}\in
C^\infty_\mathrm{c}(\Omega),$ $\delta\in(0,1),$ and $\psi\in C^\infty_\mathrm{c}\bigl(\Omega, C^\infty_\#(Y)\bigr)$ such that
\begin{equation*}
  \bigl\|g^{1,\delta}-g^1\bigr\|_{H^1(\Omega)}+\bigl\|\nabla_y\psi^{\delta}-\nabla_y\psi\bigr\|_{L^2(\Omega\times
    Y)}\leq\delta\ \ \ \ \ \forall\delta.
\end{equation*}
For $\e>0,$ $\delta\in(0,1),$ define
\begin{equation*}
  g^{1,\delta}_\e(x):=g^{1,\delta}(x)+\e\psi^{\delta}(x,x/\e),
\end{equation*}
and set
\begin{equation*}
  d_\e^\delta:=d_\e\bigl(\nabla  g^{1,\delta}_\e,\nabla
  g^1+\nabla_y\psi\bigr)+\bigl\|g^{1,\delta}_\e-g^1\bigr\|_{L^2(\Omega)}
  +c_\e\bigl\|\nabla g^{1,\delta}_\e\bigr\|_{L^\infty(\Omega)}.
\end{equation*}
By construction, we have
$\lim_{\delta\downarrow0}\limsup_{\e\downarrow 0}
d_\e^\delta=0$, and Lemma~\ref{lem:attouch} yields a function
$\e\mapsto \delta(\e)$ with $\lim_{\e\downarrow 0}d^{\delta(\e)}_\e=
0$. In view of Step~1, this implies that the diagonal sequence $\tilde g^1_\e:=g^{1,\delta(\e)}_\e$ satisfies \eqref{L4:pf0}. Now, for each $\varepsilon>0,$ let $g^1_\e$ denote the function satisfying \eqref{L1:eq:1}(ii) and such that
$g^1_\e=\tilde g^1_\e$ on $\Omega^1_\e$. To conclude the argument, we
only need to
show that $g^1_\e$ satisfies \eqref{L4:pf0}. Consider the difference $\eta_\e:=\tilde
g^1_\e-g^1_\e$. Since $\eta_\e$ is bounded in $H^1(\Omega)$ and
$\eta_\e=0$ in $\Omega^1_\e$, we have $\eta_\e\wto 0$ in
$H^1(\Omega)$, and, up to a subsequence, $\nabla\eta_\e\wtto
\nabla_y\varphi$ for some $\varphi\in L^2\bigl(\Omega,H^1_0(Y^0)\bigr)$. On
the other hand, since $g^1_\e$ satisfies \eqref{L1:eq:1} (ii) and $\tilde
g^{1}_\e$ satisfies \eqref{L4:pf0}, we deduce that
\begin{eqnarray*}
  \|\nabla\eta_\e\|_{L^2(\Omega)}^2=\int_{\Omega^0_\e}\nabla\eta_\e:\nabla\eta_\e\,dx
  =\int_{\Omega^0_\e}\nabla\tilde
  g^1_\e:\nabla\eta_\e\,dx\to\iint_{\Omega\times Y^0}(\nabla g^1+\nabla_y\psi):\nabla_y\varphi\,dxdy.
\end{eqnarray*}
Since $\nabla g^1$ is independent of $y$, and because $\psi\in
L^2\bigl(\Omega,\mathcal A(Y^0)\bigr)$, the integral on the right-hand side
vanishes. Hence, $\|\nabla\eta_\e\|_{L^2(\Omega)}^2\to 0$, and thus
$g^1_\e$ satisfies \eqref{L4:pf0}.
\smallskip

\step 3 Conclusion.

As can be shown by appealing to a combination of a density argument and a diagonal-sequence argument, similar to Step~1, there exists a
sequence $g^0_\e\in H^1_0(\Omega^0_\e)$ such that
\begin{equation*}
  g^0_\e\stto g^0,\qquad \e\nabla g^0_\e\stto \nabla_yg^0,\qquad
  \limsup\limits_{\e\downarrow 0}c_\e\bigl\|\e^{1-\gamma}\nabla g^0_\e\bigr\|_{L^\infty(\Omega)}=0.
\end{equation*}
Now define $\varphi_\e(x):=\e^{1-\gamma}g^0_\e+g^1_\e$, and note that
$(g^0_\e,g^1_\e)$ satisfy \eqref{L1:eq:1}. In view of the convergence of
$g^0_\e$ and $g^1_\e,$ the sequence $\varphi_\e$ has the required properties.\end{proof}

\subsection{Proof  of Theorems~\ref{T1}, \ref{T1b} and Proposition~\ref{P1}: small strain regime}\label{S:P2}
As a preliminary remark, we note that two different effects play a role when
passing to the limit $\e\downarrow 0$ in the small strain regime:
\begin{itemize}
\item The non-convex energy functional is linearised at identity map
  (which is a stress-free state for $\mathcal E_\e$) --  this corresponds to the passage from nonlinear to
  linearised elasticity.
\item The obtained linearised, still oscillating, convex-quadratic energy is homogenised.
\end{itemize}
The following lemma is used to treat both effects simultaneously. Its
proof combines convex homogenisation methods (e.g.~\cite[Proposition 1.3]{Visintin-07}) with a ``careful Taylor expansion''
in the spirit of \cite[Proof of Theorem~6.2]{FJM2002}.
For notational convenience, we introduce two ``linearised'' functionals: 
\begin{itemize}
\item  For $G_\e=(G^0_\e,G^1_\e)\in L^2(\Omega,\R^{d\times
    d})\times L^2(\Omega,\R^{d\times d})$ set
  \begin{equation*}
    \mathcal Q_\e(G_\e)=\mathcal
    Q_\e(G^0_\e,G^1_\e):=\int_{\Omega^0_\e}Q^0\bigl(G^0_\e(x)\bigr)\,dx+\int_{\Omega^1_\e}Q^1\bigl(G^1_\e(x)\bigr)\,dx.
  \end{equation*}
\item  For $G=(G^0,G^1)\in L^2(\Omega\times Y,\R^{d\times
    d})\times L^2(\Omega\times Y,\R^{d\times d})$ set
  \begin{equation*}
    \mathcal Q(G)={\mathcal Q}(G^0,G^1):=\int_{\Omega\times
      Y^0}Q^0\bigl(G^0(x,y)\bigr)\,dxdy+\int_{\Omega\times Y^1}Q^1\bigl(G^1(x,y)\bigr)\,dxdy.
  \end{equation*}
\end{itemize}
\begin{lemma}\label{AL3}
  Consider sequences $g^0_\e, g^1_\e\in H^1(\Omega)$ that satisfy
  \begin{equation*}
    \limsup\limits_{\e\downarrow 0}\Bigl(\bigl\|\e\nabla g^0_\e\bigr\|_{L^2(\Omega)}+\bigl\|\nabla g^1_\e\bigr\|_{L^2(\Omega)}\Bigr)<\infty.
  \end{equation*}
  Set $\varphi_\e:=\e^{1-\gamma}g^0_\e+g^1_\e$ and $G_\e=(G_\e^0, G_\e^1):=\bigl(\e\nabla
  g^0_\e,\nabla g^1_\e\bigr)$.
  \begin{enumerate}[(a)]
  \item If $G_\e\wtto G$ and $m_\e=o(\e^{\gamma})$ as $\e\downarrow0,$ then
    \begin{equation}
      \label{AL3:2a}
      \begin{aligned}
        \liminf\limits_{\e\downarrow 0}m_\e^{-2}\mathcal E_\e(m_\e\varphi_\e)\geq\liminf\limits_{\e\downarrow
          0}\mathcal Q_\e(\theta_\e G_\e)\geq \mathcal Q(G),
     \end{aligned}
    \end{equation}
    where $\theta_\e:\Omega\to\{0,1\}$ is defined by
    \begin{equation*}
      \theta_\e(x):=
      \begin{cases}
        1,&\ \mathrm{if\ } |\nabla\varphi_\e|\leq {(m_\e\e^{\gamma})}^{-1/2},\\
        0,&\ \mathrm{otherwise.}
      \end{cases}
    \end{equation*}
  \item  If $G_\e\stto G$, $m_\e=o(\e^\gamma)$ as $\e\downarrow0,$ and
    \begin{equation}
      \label{AL3:ass2}
      \limsup\limits_{\e\downarrow 0}\|m_\e\nabla\varphi_\e\|_{L^\infty(\Omega)}=0,
    \end{equation}
    then
    \begin{equation*}
      \begin{aligned}
        \lim\limits_{\e\downarrow 0}m_\e^{-2}\mathcal
        E_\e(m_\e\varphi_\e)=\mathcal Q(G).
      \end{aligned}
    \end{equation*}
  \end{enumerate}
\end{lemma}
\begin{remark}\label{R:theta}
  In Lemma~\ref{AL3}, following \cite{FJM2002}, the function $\theta_\e$ is introduced, in order to
  truncate the peaks of $G_\e$. This is needed for exploiting the
  quadratic expansion (W3). Since both $\e\nabla
  g^0_\e$ and $\nabla g^1_\e$ are assumed to be bounded sequences in
  $L^2(\Omega)$, we deduce, from the definition of $\theta_\e$, the fact that $m_\e\e^{-\gamma}=o(1)$ and
  the Chebyshev inequality, that
  \begin{equation}\label{eq:95}
    \forall r<\infty\qquad \|\theta_\e-1\|_{L^r(\Omega)}\to
    0,\qquad\mbox{and}\qquad\|\theta_\e\|_{L^\infty(\Omega)}\leq 1.
  \end{equation}
\end{remark}
In the proof of Lemma~\ref{AL3} we need to pass to the limit
in products of the form $f_\e\theta_\e$, where $\theta_\e$ satisfies
\eqref{eq:95}, or $f_\e\chi_\e$, where $\chi_\e$ denotes the indicator
of $\Omega^0_\e$. This is done by appealing to the
next two lemmas, the proofs of which are elementary and left to the reader.
\begin{lemma}
  \label{L:theta}
  Let $f_\e,\theta_\e$ be sequences in $L^2(\Omega)$ and assume that $\theta_\e$ satisfies \eqref{eq:95}, then the following implications
  are valid:
  \begin{align*}
    \limsup\limits_{\e\downarrow
      0}\|f_\e\|_{L^2(\Omega)}<\infty\qquad&\Rightarrow\qquad \|\theta_\e f_\e-f_\e\|_{L^p(\Omega)}\to 0\ \ \ \forall p<2,\\
    f_\e\wto f\qquad&\Rightarrow\qquad\theta_\e f_\e\wto f,\\
    f_\e\wtto f\qquad&\Rightarrow\qquad\theta_\e f_\e\wtto f,\\
    \left\{\begin{aligned}
      &f_\e\wto f\\
      &\bigl(|f_\e|^2\bigr)\ \text{\rm equi-integrable}    
      \end{aligned}\right\}\qquad&\Rightarrow\qquad\|\theta_\e f_\e-f_\e\|_{L^2(\Omega)}\to 0.
  \end{align*}
\end{lemma}
\begin{lemma}
  \label{L:chi}
  Suppose that $f_\e$ be a sequence in $L^2(\Omega)$ and, as above, let $\chi_\e$ denote
  the set indicator function of $\Omega^0_\e$. Then the following
  implications hold:
  \begin{align*}
    f_\e\wtto f\qquad&\Rightarrow\qquad\chi_\e f_\e\wtto \chi(y)f(x,y),\\
    f_\e\stto f\qquad&\Rightarrow\qquad\chi_\e f_\e\stto \chi(y)f(x,y),\\
  \end{align*}
  where $\chi$ denotes the indicator function of $Y^0$.  
\end{lemma}

\begin{proof}[Proof of Lemma~\ref{AL3}]
  \step 1 Linearisation.
  
  We claim that the following statement holds for $i=1,2$:
  Let $F_\e$ denote a sequence in $L^2(\Omega,\R^{d\times d})$, and
  let $c_\e$ be a sequence of positive numbers converging to zero, such that
  \begin{equation}\label{AL3:linass}
    \limsup_{\e\to 0}c_\e||F_\e||_{L^\infty(\Omega)}=0.
  \end{equation}
  Then the convergence
  \begin{equation}
    \label{AL3:pflin}
    \lim\limits_{\e\downarrow 0}\left|c_\e^{-2}\int_\Omega W^i(I+c_\e
      F_\e)\,dx-\int_{\Omega}Q^i\bigl(F_\e(x)\bigr)\,dx\right|=0
  \end{equation}
holds.  Indeed, thanks to (W3) we have
  \begin{equation*}
    \left|c_\e^{-2}W^i(I+c_\e F_\e)-Q^i(F_\e)\right|\leq
    |F_\e|^2r^i(c_\e|F_\e|)\leq |F_\e|^2r^i\bigl(c_\e\|F_\e\|_{L^\infty(\Omega)}\bigr)\qquad\mbox{a.e.}
  \end{equation*}
  Thanks to \eqref{AL3:linass}, and since $F_\e$ is bounded in
  $L^2(\Omega)$, the right-hand side converges to zero in
  $L^1(\Omega)$, and \eqref{AL3:pflin} follows.
  \smallskip

  \step 2 Proof of part (a).

  Since the energy densities $W^0,$ $W^1$ are minimised at the identity, {\it cf.} (W2), we have
  \begin{equation}\label{AL3:pf1}
   m_\e^{-2}\mathcal
    E_\e(\varphi_\e)\geq(m_\e\e^{-\gamma})^{-2}\int_\Omega
    W^0\bigl(I+m_\e\e^{-\gamma}\theta_\e F^0_\e(x)\bigr)\,dx+m_\e^{-2}\int_\Omega W^1\bigl(I+m_\e\theta_\e F^1_\e(x)\bigr)\,dx,
  \end{equation}
  where
  \begin{equation}\label{AL3:pf2}
    F^0_\e:=\chi_\e\bigl(G^0_\e+\e^\gamma G^1_\e\bigr),\qquad F^1_\e(x):=(1-\chi_\e)G^1_\e.
  \end{equation}
  Thanks to the definition of $\theta_\e$ we have
  $\bigl\|m_\e\e^{-\gamma}\theta_\e
  F^0_\e\bigr\|_{L^\infty(\Omega)}+\bigl\|m_\e\theta_\e
  F^1_\e\bigr\|_{L^\infty(\Omega)}\to 0$, so that we may apply
  \eqref{AL3:pflin} to the right-hand side in \eqref{AL3:pf1}. We get
  \begin{eqnarray*}
    \liminf\limits_{\e\downarrow 0}m_\e^{-2}\mathcal
    E_\e(\varphi_\e)\geq\liminf\limits_{\e\downarrow 0}\mathcal
    Q_\e\bigl(\theta_\e F^0_\e,\theta_\e F^1_\e\bigr)=\liminf\limits_{\e\downarrow 0}\mathcal
    Q_\e\bigl(\theta_\e G_\e^1,\theta_\e G^0_\e\bigr),
  \end{eqnarray*}
where for the last identity we used the facts that  $F^1_\e=G^1_\e$ on
$\Omega^1_\e$ and
\[
\bigl\|F^0_\e-G^0_\e\bigr\|_{L^2(\Omega^0_\e)}=\bigl\|\e^\gamma\nabla g^1_\e\bigr\|_{L^2(\Omega^0_\e)}\to
0.
\] 
It remains to argue that
\begin{equation*}
  \liminf\limits_{\e\downarrow 0}\mathcal
    Q_\e(\theta_\e G_\e)\geq
    \mathcal Q(G).
\end{equation*}
In order to show this, notice that
\begin{equation}\label{AL3:pf:5}
  Q_\e(\theta_\e G_\e)=\int_{\Omega}Q^0\bigl(\theta_\e\chi_\e G^0_\e\bigr)\,dx+\int_{\Omega}Q^1\bigl(\theta_\e(1-\chi_\e)G^1_\e\bigr)\,dx.
\end{equation}
From $G_\e\wtto G$ we deduce, using Lemma~\ref{L:theta},
Remark~\ref{R:theta} and Lemma~\ref{L:chi}, that
\begin{equation*}
  \theta_\e\chi_\e G^0_\e\wtto \chi(y) G^0(x,y),\qquad \theta_\e(1-\chi_\e)G^1_\e\wtto\bigl(1-\chi(y)\bigr)G^1(x,y).
\end{equation*}
By appealing to the lower semicontinuity of convex integral
functionals with respect to weak two-scale convergence ({\it cf.}
\cite[Proposition~1.3]{Visintin-07}), we deduce that the $\liminf$ of the
right-hand side in \eqref{AL3:pf:5} is bounded below by $\mathcal
Q(G)$. This completes the argument.
\smallskip

  \step 3 Proof of part (b).

  We claim that
  \begin{equation}\label{AL3:pf0}
    \lim\limits_{\e\downarrow 0}\left|\frac{1}{m_\e^2}\int_{\Omega^0_\e}W(\mathcal
      E_\e(m_\e\varphi_\e)-\mathcal Q_\e(\e\nabla g^0_\e,\nabla g^1_\e)\right|=0.
  \end{equation}
  Note that
  \begin{equation*}
  m_\e^{-2}\mathcal
    E_\e(\varphi_\e)=(m_\e\e^{-\gamma})^{-2}\int_\Omega
    W^0\bigl(I+m_\e\e^{-\gamma} F^0_\e(x)\bigr)\,dx + m_\e^{-2}\int_\Omega W^1\bigl(I+m_\e F^1_\e(x)\bigr)\,dx
  \end{equation*}
  where $F^0_\e$ and $F^1_\e$ are defined in \eqref{AL3:pf2}. By
  \eqref{AL3:ass2} we have \mbox{$\bigl\|m_\e\e^{-\gamma} F^0_\e\bigr\|_{L^\infty(\Omega)}+\bigl\|m_\e F^1_\e\bigr\|_{L^\infty(\Omega)}\to0$}
  and \eqref{AL3:pflin} yields
  \begin{equation*}
    \Bigl|m_\e^{-2}\mathcal E_\e(m_\e\varphi_\e)-\mathcal
      Q_\e\bigl(F^0_\e,F^1_\e\bigr)\Bigr|\to 0\ \ \mathrm{as}\ \e\downarrow0.
  \end{equation*}
  Since $G_\e\stto G$ we have, thanks to Lemma~\ref{L:chi}:
\begin{equation*}
  \chi_\e F^0_\e\stto \chi(y) G^0(x,y),\qquad (1-\chi_\e)F^1_\e\wtto\bigl(1-\chi(y)\bigr)G^1(x,y).
\end{equation*}
Hence, the continuity of convex integral
functionals with respect to strong two-scale convergence ({\it cf.} \cite{Visintin}) yields
\begin{equation*}
  \mathcal Q_\e\bigl(F^0_\e,F^1_\e\bigr)\to \mathcal Q(G),
\end{equation*}
which completes the argument.
\end{proof}
We are now in a position to prove the $\Gamma$-convergence statement
for the energies $\mathcal E_\e$.
\begin{proof}[Proof of Theorem~\ref{T1}]

\step 1 Part (a) (Compactness).

Thanks to (W2) we have 
\[
m_\e^{-2}\mathcal
E_\e(m_\e\varphi_\e)\geq c_0m_\e^{-2}\Phi^{\gamma}_{\e}(m_\e\varphi_\e).
\]
Hence, the claim of Theorem~\ref{T1}(a) directly follows from 
Lemma~\ref{L2}.
\smallskip

\step 2 Part (b) (Lower bound).

Without loss of generality we assume that
\begin{equation*}
  \liminf\limits_{\e\downarrow 0}m_\e^{-2}\mathcal
      E_\e(m_\e\varphi_\e)=\limsup\limits_{\e\downarrow 0}m_\e^{-2}\mathcal
      E_\e(m_\e\varphi_\e)<\infty.
\end{equation*}
Furthermore, thanks to Lemma~\ref{L2}, we can assume in addition that $\nabla g^1_\e\wtto \nabla g^1+\nabla_y\psi$
for some $\psi\in L^2\bigl(\Omega,\mathcal A(Y^0)\bigr)$, so that
\begin{equation*}
  G_\e:=\bigl(\e\nabla g^0_\e,\nabla g^1_\e\bigr)\wtto\bigl(\nabla_y g^0,\nabla g^1+\nabla_y\psi\bigr)=:G.
\end{equation*}
Applying Lemma~\ref{AL3}(a) yields
\begin{equation*}
  \liminf\limits_{\e\downarrow 0}m_\e^{-2}\mathcal
      E_\e(m_\e\varphi_\e)\geq \int_{\Omega\times Y^0}
  Q^0\bigl(\nabla_yg^0(x,y)\bigr)\,dxdy+\int_{\Omega\times  Y^1}Q^1\bigl(\nabla g^1(x)+\nabla_y\psi(x,y)\bigr)\,dxdy.
\end{equation*}
This completes the argument, since the right-hand side is bounded from
below by $\mathcal E_{\esmall}(g^0,g^1)$.
\smallskip

\step 3 Part (c) (Upper bound).

Choose $\psi\in L^2\bigl(\Omega,\mathcal A(Y^0)\bigr)$ such that
\begin{equation}\label{T1:pf4}
  \int_{\Omega\times Y^1}Q\bigl(\nabla
  g^1(x)+\nabla_y\psi(x,y)\bigr)\,dxdy=\int_\Omega Q^1_{\hom}\bigl(\nabla g^1(x)\bigr)\,dx.
\end{equation}
Let $\varphi_\e$ denote the sequence associated with $g^0,g^1$ and
$\psi$ via Lemma~\ref{L4} with $c_\e:=m_\e$. In view of \eqref{L4:1}, applying Lemma~\ref{AL3}~(b) yields
\begin{equation*}
  \lim\limits_{\e\downarrow 0}m_\e^{-2}\mathcal
  E_\e(m_\e\varphi_\e)=\int_{\Omega\times Y}Q^0\bigl(\nabla_y
  g^0(x,y)\bigr)\,dxdy+\int_{\Omega\times  Y^1}Q^1\bigl(\nabla
  g^1(x)+\nabla_y\psi(x,y)\bigr)\,dxdy.
\end{equation*}
It follows from \eqref{T1:pf4} that the right-hand side equals $\mathcal
E_{\esmall}(g^0,g^1)$. 
\end{proof}

\begin{proof}[Proof of Proposition~\ref{P1}]
  \step 1 {\it A priori} estimate.

  We claim that for every sequence $\varphi_\e$ in $H^1_0(\Omega)$ the
  following implication holds:
  \begin{equation}\label{P1:pf1}
    \limsup\limits_{\e\downarrow 0}I_\e(\varphi_\e)<\infty\qquad\Rightarrow\qquad
    \limsup\limits_{\e\downarrow 0}m_\e^{-2}\mathcal E_\e(m_\e\varphi_\e)<\infty.  
  \end{equation}
  Indeed, we have
  \begin{eqnarray*}
    \left|\int_\Omega \ell_\e\cdot m_\e\varphi_\e\,dx\right|
    &\leq&
    m_\e\Bigl(\|\ell_\e\|_{L^2(\Omega)}\bigl\|g_\e^1\bigr\|_{L^2(\Omega)}+\e^{1-\gamma}\|\ell_\e\|_{L^2(\Omega^0_\e)}\bigl\|g^0_\e\bigr\|_{L^2(\Omega^0_\e)}\Bigr)\\
    &\stackrel{\eqref{eq:79}}{\leq}&
    m_\e\Bigl(\bigl\|m_\e g_\e^1\bigr\|_{L^2(\Omega)}+\bigl\|m_\e
      g^0_\e\bigr\|_{L^2(\Omega)}\Bigr)\,\stackrel{\eqref{L1:eq:2}}{\lesssim}\,
    m_\e^2\sqrt{m_\e^{-2}\Phi^\gamma_\e(m_\e\varphi_\e)}\\
    &\stackrel{\text{(W2)}}{\lesssim}&
    m_\e^2\sqrt{m_\e^{-2}\mathcal E_\e(m_\e\varphi_\e)}.
  \end{eqnarray*}
  Combining this with the definition
  of $I_\e$ we get
  \begin{equation*}
  m_\e^{-2}\mathcal E_\e(m_\e\varphi_\e)\lesssim
    I_\e(\varphi_\e)+\sqrt{m_\e^{-2}\mathcal E_\e(\varphi_\e)},
  \end{equation*}
  which implies \eqref{P1:pf1}. \smallskip

  \step 2 The proof of parts (a) and (b).

  The existence of a minimiser to $I_{\esmall}$ follows by the direct
  method. The minimiser $(g^0_*,g^1_*)$ is unique, since the implication
  \begin{equation*}
    \int_{\Omega\times Y}Q^0\bigl(\nabla_y\tilde g^0(x,y)\bigr)\,dxdy+\int_\Omega
    Q^1_{\hom}\bigl(\nabla\tilde g^1(x)\bigr)\,dx=0\qquad\Rightarrow\qquad \tilde g^0=\tilde g^1=0
  \end{equation*}
  holds for all $\tilde g^0\in L^2\bigl(\Omega,H^1_0(Y^0)\bigr)$ and $\tilde
  g^1\in H^1_0(\Omega)$.

  The remaining claims of Proposition~\ref{P1} follow from the
  standard $\Gamma$-convergence arguments ({\it cf.} \cite[Corollary 7.20]{DalMaso}), provided the  functionals $I_\e$, $\e>0$, are equi-coercive and $\Gamma$-converge to
  $I_{\esmall}$. Indeed, thanks to \eqref{eq:90}, it is easy to check
  that $\varphi_\e\wtto (g^0,g^1)$ implies
  \begin{equation}\label{P1:pf0}
    \lim\limits_{\e\downarrow 0}\frac{1}{m_\e^2}\int_\Omega \ell_\e\cdot
    m_\e\varphi_\e\,dx= \int_\Omega\biggl(\int_{Y^0}\ell^0\cdot
      g^0\,dy+\ell^1\cdot g^1\biggr)\,dx,
  \end{equation}
  since the integral on the left-hand side only involves products of  weakly and strongly two-scale
  convergent factors, {\it cf.} \cite[Proposition 2.8]{Visintin}).
  In combination with Theorem~\ref{T1}, this implies that $I_\e$
  $\Gamma$-converges to $I_{\esmall}$. In addition, the trivial
  inequality
  \begin{equation*}
    \inf_{\varphi\in H^1_0(\Omega)}I_\e(\varphi)\leq I_\e(0)=0,
  \end{equation*}
  combined with \eqref{P1:pf1} and Lemma~\ref{L2}, proves that the
  functionals $I_\e$ are equi-coercive.
\end{proof}

For the proof of Theorem~\ref{T1b} we make use of the following lemma:
\begin{lemma}[Decomposition Lemma, see \cite{Fonseca},
  \cite{Kristensen}]\label{AL:decomp}
  Let $v\in H^1_{0}(\Omega)$ and let $v_\e\in H^1(\Omega)$ be a sequence with $v_\e\wto
  v$ in $H^1(\Omega)$. Then there exists a sequence
  $\overline V_\e\in H^1(\Omega)$ such that the following
  properties hold for a subsequence of $v_\varepsilon$ (not relabelled):
  \begin{enumerate}[(a)]
  \item $\overline V_\e\wto v$ in $H^1(\Omega);$
  \item $\overline V_\e=v_\e$ in a neighbourhood of $\partial\Omega;$
  \item $\bigl(|\nabla \overline V_\e|^2\bigr)$ is equi-integrable;
  \item $\bigl|\{\,x\in\Omega\,:\,v_\e(x)\neq \overline V_\e(x)\,\}\bigr|\to 0$ as $\e\downarrow0.$
  \end{enumerate}
\end{lemma}

\begin{proof}[Proof of Theorem~\ref{T1b}]
It suffices to prove the theorem for a subsequence. Throughout the proof we write
\begin{equation*}
  G_\e:=\bigl(\e\nabla g^0_\e,\nabla g^1_\e\bigr),\qquad G:=\bigl(\nabla_y
  g^0_*,\nabla g^1_*+\nabla_y\psi_*\bigr).
\end{equation*}
Furthermore, we make use of the functionals $\mathcal Q_\e$ and
$\mathcal Q$ introduced at the beginning of Section~\ref{S:P2}. Recall
that
\begin{equation*}
  \mathcal E_{\esmall}\bigl(g^0_*,g^1_*\bigr)=Q(G).
\end{equation*}

\step 1 Convergence of $\varphi_\e$ and of the corresponding energy values.

We claim that, as $\e\downarrow0,$ one has
\begin{align}
  \label{T1b:pf5a}
  &G_\e\wtto G,\\
  \label{T1b:pf5c}
  &I_\e(\varphi_\e)\to I_{\esmall}\bigl(g^0_*,g^1_*\bigr),\\[0.35em]
  \label{T1b:pf5d}
  &m_\e^{-2}\mathcal E_\e(m_\e\varphi_\e)\to \mathcal E_{\esmall}\bigl(g^0_*,g^1_*\bigr).
\end{align}
Indeed, from
Proposition~\ref{P1} we immediately deduce that $\varphi_\e\wtto
\bigl(g^0_*,g^1_*\bigr)$ and \eqref{T1b:pf5c}. Furthermore, in view of the continuity of the loading term, {\it cf.}
\eqref{P1:pf0}, this implies \eqref{T1b:pf5d}. For \eqref{T1b:pf5a},
it remains to argue that $\nabla g^1_\e\wtto \nabla
g^1_*+\nabla_y\psi_*$.  Thanks to $\varphi_\e\wtto (g^0_*,g^1_*)$ and Lemma~\ref{L2} we have, up to a subsequence, $\nabla
g^1_\e\wtto \nabla g^1_*+\nabla_y\psi$ for some $\psi\in
L^2(\Omega,\mathcal A(Y^0))$. Furthermore, from \eqref{T1b:pf5d} and Lemma~\ref{AL3}~(a)
we infer that
\begin{equation*}
  E_{\esmall}(g^*_0,g^*_1)=\liminf\limits_{\e\downarrow 0}m_\e^{-2}\mathcal
  E_\e(m_\e\varphi_\e)\geq \int_{\Omega\times Y^0}Q^0\bigl(\nabla_y
  g^0_*\bigr)\,dxdy+\int_{\Omega\times Y^1}Q^1\bigl(\nabla g^1_*+\nabla_y\psi\bigr)\,dxdy.
\end{equation*}
This, in particular, implies
\begin{equation*}
  \int_{\Omega\times Y^1}Q^1\bigl(\nabla
  g^1_*+\nabla_y\psi\bigr)\,dxdy=\int_\Omega Q^1_{\hom}\bigl(\nabla g^1_*\bigr)\,dx.
\end{equation*}
In view of \eqref{eq:91} we conclude that $\psi=\psi_*$ and \eqref{T1b:pf5a} follows.
\smallskip

\step 2 Equi-integrable decomposition.

We claim that for a subsequence
(not relabelled) there exist sequences $\bar
g^0_\e,\overline G^1_\e\in H^1_0(\Omega)$ such that $\overline G_\e:=(\e\nabla\overline G^0_\e,\nabla\overline G^1_\e)$ satisfies
\begin{equation}\label{T1b:pf2-1}
  \overline G_\e-G_\e\wtto 0,\qquad ||\overline G_\e-G_\e||_{L^p(\Omega)}\to
  0,
\end{equation}
and
\begin{equation}
  \label{T1b:pf2-2}
  \mathcal Q_\e(\overline G_\e)\to\mathcal E_{\esmall}\bigl(g^0_*,g^1_*\bigr).
\end{equation}
To show the above, notice that thanks to Lemma~\ref{AL:decomp} there exist
sequences $\overline G^0_\e, \overline G^1_\e\in H^1_0(\Omega)$ such that
\begin{itemize}
\item $\bigl(\e^2|\nabla\overline G^0_\e|^2\bigr)$ and $\bigl(|\nabla\bar
g^1_\e|^2\bigr)$ are equi-integrable,
\item the indicator function $\bar\theta_\e$ defined by
  \begin{equation}\label{T1b:pf7}
    \bar\theta_\e(x):=
    \begin{cases}
      1,&\text{if}\,\,\,\overline G^0_\e(x)=g^0_\e(x)\mbox{ and }\overline G^1_\e(x)=g^1_\e(x),\\
      0,&\text{otherwise},
    \end{cases}
  \end{equation}
  satisfies \eqref{eq:95}.
\end{itemize}
Since $\overline G_\e-G_\e=(1-\bar\theta_\e)(\overline G_\e-G_\e)$ (and $p<2$),
the convergence \eqref{T1b:pf2-1} follows from the boundedness of the sequence $(\overline G_\e-G_\e)$ in $L^2(\Omega)$, Lemma~\ref{L:theta}, and H\"older's inequality.

We prove \eqref{T1b:pf2-2}. Thanks to \eqref{T1b:pf5a} 
we have $\overline G_\e\wtto G$, so that (due to the lower semicontinuity of convex integral functionals with respect to weak two-scale convergence, {\it cf.} \cite[Proposition~1.3]{Visintin-07}):
\begin{equation*}
  \liminf\limits_{\e\downarrow 0}\mathcal Q_\e(\overline G_\e)\geq \mathcal
  Q(G)=\mathcal E_{\esmall}\bigl(g^0_*,g^1_*\bigr).
\end{equation*}
Hence, for \eqref{T1b:pf2-2} it suffices to prove the opposite
estimate, i.e.~$\limsup_{\e\downarrow 0}\mathcal Q_\e(\overline G_\e)\leq\mathcal
E_{\esmall}(G)$, which, thanks to \eqref{AL3:2a} and \eqref{T1b:pf5d}, follows from
\begin{equation}\label{T1b:pf6}
  \limsup\limits_{\e\downarrow 0}\left(\mathcal Q_\e(\overline G_\e)-\mathcal Q_\e(\theta_\e
    G_\e)\right)\leq 0.
\end{equation}
In order to show \eqref{T1b:pf6} notice that since the supports of
$\bar\theta_\e$ and $(1-\bar\theta_\e)$ are disjoint, and because
$\bar\theta_\e\theta_\e G_\e=\bar\theta_\e\theta_\e\overline G_\e$ (cf. \eqref{T1b:pf7}), an expansion of the squares yields
\begin{eqnarray*}
  \mathcal Q_\e(\overline G_\e)-\mathcal Q_\e(\theta_\e G_\e)
  &=&
  \mathcal Q_\e(\bar\theta_\e\overline G_\e)+\mathcal
  Q_\e\bigl((1-\bar\theta_\e)\overline G_\e\bigr)
  -\mathcal Q_\e(\bar\theta_\e\theta_\e G_\e)-\mathcal
  Q_\e\bigl((1-\bar\theta_\e)\theta_\e G_\e\bigr) \\
  &=&
  \mathcal Q_\e\bigl(\bar\theta_\e(1-\theta_\e)\overline G_\e\bigr)+\mathcal
  Q_\e\bigl((1-\bar\theta_\e)\overline G_\e\bigr)-\mathcal
  Q_\e\bigl((1-\bar\theta_\e)\theta_\e G_\e\bigr)\\
  &\leq&
  \mathcal Q_\e\bigl(\bar\theta_\e(1-\theta_\e)\overline G_\e\bigr)+\mathcal
  Q_\e\bigl((1-\bar\theta_\e)\overline G_\e\bigr).
\end{eqnarray*}
It is easy to check that $\bar\theta_\e(1-\theta_\e)$ and
$1-\bar\theta_\e$ converge to zero in $L^r(\Omega)$ for all $r<\infty$. Hence, since $|\overline G_\e|^2$  is equi-integrable,
Lemma~\ref{L:theta} implies that the right-hand side of the previous
estimate converges to zero, and \eqref{T1b:pf6} follows.
\smallskip

\step 3 Error estimate.

We claim that
\begin{equation}\label{T1b:pf9-1}
  \int_{\Omega^0_\e}\e\Bigl|\sym\nabla\bigl(\bar
  g^0_\e-g^0_{*,\e}\bigr)\Bigr|^2\,dx+\int_{\Omega^1_\e}\Bigl|\sym\nabla\bigl(\bar
  g^1_\e-g^1_{*,\e}\bigr)\Bigr|^2\,dx\to 0\ \ \ \mathrm{as}\ \e\downarrow0.
\end{equation}
For the argument set $G_{*,\e}:=\bigl(\e\nabla g^0_{*,\e},\nabla
g^1_{*,\e}\bigr)$. In view of \eqref{eq:sym} it suffices to argue that
\begin{equation*}
  \mathcal Q_\e(\overline G_\e-G_{*,\e})\to 0\ \ \ \mathrm{as}\ \e\downarrow0.
\end{equation*}
The latter can be seen as follows: We have
\begin{eqnarray*}
  \mathcal Q_\e(\overline G_\e-G_{*,\e})=\mathcal Q_\e(\overline G_\e)-\mathcal
  Q_\e(G_{*,\e})+2\mathcal B_\e(G_{*,\e};G_{*,\e}-\overline G_\e),
\end{eqnarray*}
where $\mathcal B_\e$ denotes the bilinear form associated with
$\mathcal Q_\e$.

The difference of the two quadratic terms on the right-hand side
converges to zero, since $G_{*,\e}$ is associated with a recovery
sequence, and thanks to \eqref{T1b:pf2-2}. On the other hand, since $G_{*,\e}$ strongly two-scale converges, and $G_{*,\e}-\bar
G_\e\wtto 0$ by \eqref{T1b:pf2-1}, we deduce that
\begin{equation*}
  \lim\limits_{\e\downarrow 0}\mathcal B_\e(G_{*,\e};G_{*,\e}-\overline G_\e)=0,
\end{equation*}
as $\mathcal B_\e(G_{*,\e};G_{*,\e}-\overline G_\e)$ only involves products
between a weakly and a strongly two-scale convergent factor ({\it cf.}
\cite[Proposition 2.8]{Visintin}).
\smallskip

\step 4 Conclusion  (Proof of \eqref{T1b:1}).

We split the estimate into
\begin{align}
  \label{L4:1a}
  &\int_{\Omega}\Bigl(\bigl|g^0_\e-g^0_{*,\e}\bigr|^p+\bigl|\e\nabla g^0_\e-\e\nabla
  g^0_{*,\e}\bigr|^p\Bigr)\,dx\to 0\ \ \ \mathrm{as}\ \e\downarrow0,\\
  \label{L4:1b}
  &\int_{\Omega}\Bigl(\bigl|g^1_\e-g^1_{*,\e}\bigr|^p+\bigl|\nabla g^1_\e-\nabla
  g^1_{*,\e}\bigr|^p\Bigr)\,dx\to 0\ \ \ \mathrm{as}\ \e\downarrow0.
\end{align}
Thanks to \eqref{T1b:pf2-1} and Step~3 we have
\begin{equation}\label{T1b:pf8-1}
  \int_{\Omega^0_\e}\e\Bigl|\sym\nabla\bigl(g^0_\e-g^0_{*,\e}\bigr)\Bigr|^p\,dx+\int_{\Omega^1_\e}\Bigl|\sym\nabla\bigl(g^1_\e-g^1_{*,\e}\bigr)\Bigr|^p\,dx\to 0 \ \ \ \mathrm{as}\ \e\downarrow0.
\end{equation}
Argument for \eqref{L4:1a}: Set $\eta_\e^0:=g^0_\e-g^0_{*,\e}$. Since
$\eta_\e^0\in H^1_0(\Omega^0_\e)\subset H^1_0(\Omega)$, Korn's
inequality yields
\begin{equation*}
  \int_{\Omega}\bigl|\nabla \eta^0_\e\bigr|^p\,dx\leq
  C\int_{\Omega}\bigl|\sym\nabla\eta^0_\e\bigr|^p\,dx= C\int_{\Omega^0_\e}\bigl|\sym\nabla\eta^0_\e\bigr|^p\,dx,
\end{equation*}
where $C>0$ only depends on $\Omega$, $p$, and $d$. Combined with the
improved Poincar\'e inequality \eqref{improvedPoincare} and
\eqref{T1b:pf8-1}, \eqref{L4:1a} follows.

Argument for \eqref{L4:1b}: 
We claim that \eqref{L4:1b}
follows from
\begin{equation}\label{T1b:pf8-2}
\bigl\|\sym\nabla\eta^1_\e\bigr\|_{L^p(\Omega^0_\e)}\to 0,
\end{equation}
where  $\eta_\e^1:=g^1_\e-g^1_{*,\e}$. Indeed, since $\eta_\e^1$ vanishes on
$\partial\Omega$, \eqref{T1b:pf8-2}, \eqref{T1b:pf8-1} and  Korn's
first inequality yield  \eqref{L4:1b}. 

Thanks to the definition of $\Omega^0_\e$, the argument for
\eqref{T1b:pf8-2} can be reduced to the following statement: For all $\xi\in Z_\e:=\{\,\xi\in\Z^d\,:\,\e(\xi+Y)\subset\Omega\,\}$
we have
\begin{equation}\label{T1b:pf8-3}
  \int_{\e(\xi+Y^0)}|\sym\nabla\eta_\e|^p\,dx\lesssim
  \int_{\e(\xi+Y^1)}|\sym\nabla\eta_\e|^p\,dx.
\end{equation}
For the argument consider the rescaled
function
\begin{equation*}
  \hat\eta_\e:Y\mapsto \R^d,\qquad\hat\eta_\e(y):=\eta_\e\bigl(\e(\xi+y)\bigr)-Sy+c
\end{equation*}
where $S\in\R^{d\times d}_\mathrm{skew}$ and $c\in\R^d$ are chosen
such that the Poincar\'e and Korn inequalities yield
\begin{equation}\label{T1b:pf:8-4}
  \int_{Y^1}\bigl(|\hat\eta_\e|^p+|\nabla\hat\eta_\e|^p\bigr)\,dy\lesssim\int_{Y^1}|\sym\nabla\hat\eta_\e|^p\,dy.
\end{equation}
Since both $g^1_\e$ and $g^1_{*,\e}$ satisfy \eqref{L1:eq:1}(ii), we
have $-\triangle\hat\eta_\e=0$ in $Y^0$ in the distributional sense.
Hence, thanks to Assumption~\ref{ass:reg} and \eqref{T1b:pf:8-4}, we have
\begin{equation*}
  \int_{Y^0}|\nabla\hat\eta_\e|^p\,dy\lesssim  \int_{Y^1}|\sym\nabla\hat\eta_\e|^p\,dy,
\end{equation*}
and thus
\begin{eqnarray*}
  \int_{\e(\xi+Y^0)}|\sym\nabla\eta_\e|^p\,dx&=&\e^{d-p}\int_{Y^0}|\sym\nabla\hat\eta_\e|^p\,dy\lesssim\e^{d-p}\int_{Y^1}|\sym\nabla\hat\eta_\e|^p\,dy\\
  &=&\int_{\e(\xi+Y^1)}|\sym\nabla\eta_\e|^p\,dx.
\end{eqnarray*}
\end{proof}

\subsection{Proof of Theorem~\ref{T2}: finite strain regime}
\label{finitestrain}

We define, for $g^0_\e\in
H^1_0(\Omega^0_\e)$, $g^0\in L^2\bigl(\Omega,H^1_0(Y^0)\bigr)$, and $g^1_\e,g^1\in H^1_0(\Omega),$ the following functionals:
\begin{eqnarray*}
  I^0_\e(g^0_\e)&:=&\int_{\Omega^0_\e}\Bigl(W^0\bigl(I+\e\nabla
  g^0_\e(x)\bigr)-\e^{1-2\gamma}\ell_\e\cdot
  g^0_\e\Bigr)\,dx,\\
  I^0_0(g^0)&:=&\int_{\Omega\times Y^0}\Big(\mathcal QW^0\bigl(I+\nabla_y
  g^0(x,y)\bigr)-\ell^0\cdot g^0\Bigr)\,dxdy,\\
  I^1_\e(g^1_\e)&:=&\e^{-2\gamma}\int_{\Omega^1_\e}\Bigl(W^1\bigl(I+\e^\gamma\nabla
  g^1_\e(x)\bigr)-\int_\Omega\e^{-\gamma}\ell_\e\cdot g^1_\e\Bigr)\,dx,\\
  I^1_0(g^1)&:=&\int_{\Omega}\Bigl(Q_{\hom}^1\bigl(\nabla g^1(x)\bigr)-\ell^1\cdot g^1\Bigr)\,dx.
\end{eqnarray*}
Thanks to the Lipschitz condition \eqref{eq:assW0b}, we can decompose $I_\e$ into the sum $I^0_\e+I^1_\e$ at the expense of a small error. More precisely, the following lemma holds.
\begin{lemma}
  \label{AL5} Suppose that $m_\e=\e^\gamma$.  There exists a constant $C>0$ such that for all $\e>0$ and $\varphi_\e\in H^1_0(\Omega)$,
  $\e^{1-\gamma}g^0_\e+g^1_\e\decdef\varphi_\e$ we have
  \begin{equation*}
    \Bigl|I_\e(\varphi_\e)-\bigl(I_\e^0(g^0_\e)+I^1_\e(g^1_\e)\bigr)\Bigr|\leq C
    \e^\gamma\bigl(1+\Phi^\gamma_\e(\varphi_\e)\bigr).
  \end{equation*}
 \end{lemma}
 \begin{proof}
   Note that
   \begin{align*}
     &\Bigl|I_\e(\varphi_\e)-\bigl(I_\e^0(g^0_\e)+\,I^1_\e(g^1_\e)\bigr)\Bigr|\\
     &\qquad =\biggl|\int_{\Omega^0_\e}W^0\Bigl(I+\e^{\gamma}\bigl(\e^{1-\gamma}
       \nabla g^0_\e+\nabla g^1_\e\bigr)\Bigr)-W^0\bigl(I+\e\nabla
       g^0_\e\bigr)\,dx\biggr|.
   \end{align*}
   In view of \eqref{eq:assW0b} and
   \eqref{L1:eq:2}, the statement follows.
 \end{proof}
The following lemma is a simple consequence of \cite[Lemma~21,
Lemma~22]{CherCher} and \eqref{eq:90}:
\begin{lemma}\label{AL7}
  Assume \eqref{eq:79}, \eqref{eq:90} and $m_\e=e^\gamma$.
  \begin{enumerate}[(a)]
  \item Consider a sequence $g^0_\e\in H^1_0(\Omega^0_\e)$. If
    $g^0_\e\stto g^0$ and $\e\nabla g^1_0\wtto \nabla_yg^0$, then
    \begin{equation*}
      \liminf\limits_{\e\downarrow
        0}I^0_\e\bigl(g^0_\e\bigr)\geq I^0_0\bigl(g^0\bigr).
    \end{equation*}
  \item For all $g^0\in L^2\bigl(\Omega,H^1_0(Y^0)\bigr)$ there exists a
    sequence $g^0_\e\in H^1_0(\Omega^0_\e)$ such that $g^0_\e\stto
    g^0$, $\e\nabla g^0_\e\wtto\nabla_y g^0,$ and
    \begin{equation*}
      \lim\limits_{\e\downarrow
        0}I^0_\e\bigl(g^0_\e\bigr)=I^0_0\bigl(g^0\bigr).
    \end{equation*}
  \end{enumerate}
\end{lemma}
For the stiff part one can prove (similar to Lemma~\ref{AL3}) the
following lemma:
\begin{lemma}\label{AL4}
  Assume that \eqref{eq:79}--\eqref{eq:90} hold.
  \begin{enumerate}[(a)]
  \item Consider $g^1_\e\in H^1_0(\Omega)$. If $g^1_\e\wto g^1$ weakly
    in $H^1(\Omega)$, then
    \begin{equation*}
      \liminf\limits_{\e\downarrow 0}I^1_\e\bigl(g^1_\e\bigr)\geq I^1_0\bigl(g^1\bigr).
    \end{equation*}
    \item For all $g^1\in H^1_0(\Omega)$ there exists a sequence
      $g^1_\e\in H^1_0(\Omega)$ such that
      \begin{equation*}
        g^1_\e\wto g^1\mbox{ weakly in
          $H^1(\Omega)$},\qquad\mbox{and}\qquad\lim\limits_{\e\downarrow
          0}I^1_\e\bigl(g^1_\e\bigr)= I^1_0\bigl(g^1\bigr).
      \end{equation*}
  \end{enumerate}
\end{lemma}
We proceed to the proof of Theorem~\ref{T2}.
\begin{proof}[Proof of Theorem~\ref{T2}]

  \step 1 The proof of parts (a) and (b)

  Statement (a) and (b) directly follow from Lemma~\ref{AL5}, Lemma~\ref{AL7}
  and Lemma~\ref{AL4}. 

  \step 2 Proof of \eqref{T2:1}: convergence of the minima.

  For brevity set
  \begin{align*}
    e_\e:=\inf_{\varphi\in
      H^1_0(\Omega)}\e^{-2\gamma}I_\e(\varphi_\e),\qquad
    e_0:=\inf_{g^0\in L^2(\Omega,H^1_0(Y^0))\atop g^1\in H^1_0(\Omega)} I_{\efin}(g^0,g^1).
  \end{align*}
  We prove \eqref{T2:1} in the form of the two inequalities
  \begin{equation}
    \label{T2:pf1}
    \limsup\limits_{\e\downarrow 0}e_\e\leq e_0,\qquad
    \liminf\limits_{\e\downarrow 0}e_\e\geq e_0.
  \end{equation}
  The argument for the first inequality in \eqref{T2:pf1} is standard: for $\delta>0$
  choose $(g^0,g^1)$ with $I_{\efin}(g^0,g^1)\leq e_0+\delta$. By part
  (b) there exists a recovery sequence $\varphi_\e,$ so that
   $I_\e(\varphi_\e)\to
  I_{\efin}(g^0,g^1)$. Hence
  \begin{equation*}
    \limsup\limits_{\e\downarrow 0}e_\e\leq \limsup\limits_{\e\downarrow
      0}I_\e(\varphi_\e)=I_{\efin}(g^0,g^1)\leq e_0+\delta.
  \end{equation*}
  Since this is valid for all $\delta>0$, the first inequality in \eqref{T2:pf1}
  follows.

  Next, we prove the second inequality in \eqref{T2:pf1}. Let
  $\varphi_\e$ denote a sequence with the property $\liminf_{\e\downarrow 0}e_\e=\liminf_{\e\downarrow 0}I_\e(\varphi_\e)$;
  e.g.~choose $\varphi_\e\in H^1_0(\Omega)$ such that
  $I_\e(\varphi_\e)\leq e_\e+\e$. Combining this with Lemma~\ref{AL5}
  we deduce that
  \begin{equation*}
    \liminf\limits_{\e\downarrow 0}e_\e=    \liminf\limits_{\e\downarrow
      0}I_\e(\varphi_\e)=\liminf\limits_{\e\downarrow
      0}\Bigl(I^0_\e\bigl(g^0_\e\bigr)+I^1_\e\bigl(g^1_\e\bigr)\Bigr)\geq \liminf\limits_{\e\downarrow
      0}I^0_\e\bigl(g^0_\e\bigr)+\liminf\limits_{\e\downarrow
      0}I^1_\e\bigl(g^1_\e\bigr),
  \end{equation*}
  where $\e^{1-\gamma}g^0_\e+g^1_\e\decdef\varphi_\e$. By passing to a
  subsequence, we assume without loss of generality that
  $\varphi_\e\wtto (g^0,g^1)$. Since, thanks to Lemma~\ref{AL4}, we have 
  \[
  \liminf_{\e\downarrow
    0}I^1_\e\bigl(g^1_\e\bigr)\geq I^1_0\bigl(g^1\bigr)\geq \inf_{\tilde g^1\in
    H^1_0(\Omega)}I^1_0\bigl(\tilde g^1\bigr), 
    \]
    it remains to argue that
  \begin{equation}\label{T2:pf5}
    \liminf\limits_{\e\downarrow 0}I^0_\e\bigl(g^0_\e\bigr)\geq \inf_{\tilde
      g^0\in L^2(\Omega,H^1_0(Y^0))}I^0_0\bigl(\tilde g^0\bigr).
  \end{equation}
We identify $g^0_\e$ with its extension by zero to $\R^d$, and consider the \textit{periodic unfolding} of $g^0_\e$ defined as
  \begin{equation*}
    G^0_\e:\Omega\times Y\to\R^d,\qquad 
    G^0_\e(x,y):=g^0_\e\bigl(\e\lfloor x/\e\rfloor+\e y\bigr),
  \end{equation*}
  where $\lfloor z\rfloor$ stands the unique vector in $\Z^d$ such that
  $z-\lfloor z\rfloor\in[0,1)^d$. Further, note that $G^0_\e\in
  L^2\bigl(\Omega,H^1_0(Y^0)\bigr),$ and for $\xi\in
  Z_\e:=\{\,\xi\in\Z^d\,:\,\e(\xi+Y)\subset\Omega\,\}$ and $y\in Y$ one has
  \begin{equation*}
    g^0_\e\bigl(\e(\xi+y)\bigr)=G^0_\e\bigl(\e(\xi+y),y\bigr),\qquad \e \nabla g^0_\e\bigl(\e(\xi+y)\bigr)=\nabla_y G^0_\e\bigl(\e(\xi+y),y\bigr).
  \end{equation*}
  Now we consider $I^0_\e\bigl(g^0_\e\bigr)$, which involves an integral over
  the set $\Omega^0_\e$. Since the latter can be written as a union of
  sets of the form  $\e(\xi+Y^0)$ with $\xi\in
  Z_\e$, an elementary calculation shows that
  \begin{eqnarray*}
    I_\e\bigl(g^0_\e\bigr)&=&\sum_{\xi\in Z_\e}\int_{\e(\xi+Y^0)}\Bigl(W^0\bigl(I+\e\nabla
    g^0_\e(x)\bigr)-\e^{2\gamma-1}\ell_\e(x)\cdot g^0_\e(x)\Bigr)\,dx\\
    &=&
    \sum_{\xi\in Z_\e}\e^d\int_{Y^0}\biggl(W^0\Bigl(I+\nabla_y
    G^0_\e\bigl(\e(\xi+y),y\bigr)\Bigr)-\e^{2\gamma-1}\ell_\e\bigl(\e(\xi+y)\bigr)\cdot G^0_\e\bigl(\e(\xi+y),y\bigr)\biggr)\,dy\\
    &=&
    \int_\Omega\int_{Y^0}\Bigl(W^0\bigl(I+\nabla_y
    G^0_\e(x,y)\bigr)-\e^{2\gamma-1}\ell_\e\bigl(\e(\lfloor
    x/\e\rfloor+y)\bigr)\cdot G^0_\e(x,y)\Bigr)\,dydx
  \end{eqnarray*}
  Thanks to \eqref{eq:90}, the characterisation of strong two-scale
  convergence introduced in Step~1 of the proof of Lemma~\ref{L4}, and the fact that $\limsup_{\e\downarrow 0}\bigl\|g^0_\e\bigr\|_{L^2(\Omega)}<\infty$, we have
  \begin{equation*}
    \limsup\limits_{\e\downarrow 0}\left|\int_{\Omega\times Y^0}\Big(\e^{2\gamma-1}\ell_\e\bigl(\e(\lfloor
    x/\e\rfloor+y)\bigr)-\ell^0(x,y)\Big)\cdot G^0_\e\bigl(\e(\lfloor x/\e\rfloor+y),y\bigr)\,dxdy\right|=0.
  \end{equation*}
  Hence, one has
  \begin{align*}
    &\liminf\limits_{\e\downarrow
      0}I_\e\bigl(g^0_\e\bigr)\\
    &\qquad=\liminf\limits_{\e\downarrow 0}
    \int_\Omega\int_{Y^0}\Bigl(W^0\bigl(I+\nabla_y
    G^0_\e(x,y)\bigr)-\ell^0(x,y)\cdot G^0_\e(x,y)\Bigr)\,dydx\\
    &\qquad\geq \inf_{g^0\in L^2(\Omega,H^1_0(Y^0))}
    \int_\Omega\int_{Y^0}\Bigl(W^0\bigl(I+\nabla_y
    g^0(x,y)\bigr)-\ell^0(x,y)\cdot g^0(x,y)\Bigr)\,dydx\\
    &\qquad\geq \inf_{g^0\in L^2(\Omega,H^1_0(Y^0))}
    \int_\Omega\int_{Y^0}\Bigl(\mathcal QW^0\bigl(I+\nabla_y
    g^0(x,y)\bigr)-\ell^0(x,y)\cdot g^0(x,y)\Bigr)\,dydx,
  \end{align*}
which proves \eqref{T2:pf5}.  

\step 3 The proof of the convergence \eqref{T2:2}.

Since $\mathcal Q W^0$ is quasiconvex and $Q^1_{\hom}$ quadratic, there exists a
pair $(g^0_*,g^1_*)$ that minimises $I_{\efin}$. Now the sequence
associated with $(g^0_*,g^1_*)$ via Theorem~\ref{T2}~(b) satisfies \eqref{T2:2}.
\end{proof}

\section*{Acknowledgements}
M.\,C. and K.\,C. acknowledge financial support of the Engineering and Physical Sciences Research Council 
(Grants
EP/F03797X/1 ``Variational convergence for non-linear high-contrast homogenisation problems'',
 EP/H028587/1 ``Rigorous derivation of moderate and high-contrast nonlinear composite plate theories'', and EP/L018802/1 ``Mathematical foundations of metamaterials: homogenisation, dissipation and operator theory'').
 S.\,N. acknowledges support of the German Research Foundation (Excellence Initiative).


\end{document}